\documentclass[11pt,A4paper]{article}

\usepackage[left=30mm,right=30mm,top=25mm,bottom=30mm]{geometry}
\usepackage{float}
\usepackage{xfrac}
\usepackage[normalem]{ulem}
\usepackage{longtable}
\usepackage{tikz-cd}
\usepackage{subcaption}
\usetikzlibrary{calc,intersections,angles,quotes}
\usetikzlibrary{decorations.markings}
\tikzset{
 mid arrow/.style={
 decoration={
 markings,
 mark=at position 0.5 with {\arrow[scale=0.8pt]{>}}
 },
 postaction={decorate}
 }
}
\usetikzlibrary{calc}
\usepackage[colorinlistoftodos]{todonotes}

\usepackage[percent]{overpic}
\usepackage{color}
\usepackage{amsthm,amsmath,amssymb}
\usepackage{aliascnt}
\usepackage{booktabs}
\usepackage{mathpazo}
\usepackage{microtype}
\usepackage{overpic}
\usepackage{bm}
\usepackage{sectsty}
\usepackage[pdftitle={PDFTitle},pdfauthor={Hugo Parlier},ocgcolorlinks,linkcolor=linkred,citecolor=linkred,urlcolor=linkblue]{hyperref}

\usepackage[capitalize,nameinlink,noabbrev]{cleveref}	

\usepackage{ifpdf}

\usepackage{mathtools}

\definecolor{linkred}{RGB}{128,0,128}
\definecolor{linkblue}{RGB}{16, 78, 139}

\usepackage[hang,flushmargin]{footmisc}
\usepackage{enumitem}
\usepackage{titlesec}
	\titlespacing{\section}{0pt}{12pt}{0pt}
	\titlespacing{\subsection}{0pt}{6pt}{0pt}

\titlelabel{\thetitle.\quad}

\makeatletter 

\long\def\@footnotetext#1{%
\H@@footnotetext{%
\ifHy@nesting 
\hyper@@anchor{\@currentHref}{#1}%
\else 
\Hy@raisedlink{\hyper@@anchor{\@currentHref}{\relax}}#1%
\fi 
}}

\def\@footnotemark{%
\leavevmode 
\ifhmode\edef\@x@sf{\the\spacefactor}\nobreak\fi 
\H@refstepcounter{Hfootnote}%
\hyper@makecurrent{Hfootnote}%
\hyper@linkstart{link}{\@currentHref}%
\@makefnmark 
\hyper@linkend 
\ifhmode\spacefactor\@x@sf\fi 
\relax 
}%

\ifFN@multiplefootnote%
\renewcommand*\@footnotemark{%
\leavevmode 
\ifhmode 
\edef\@x@sf{\the\spacefactor}%
\FN@mf@check 
\nobreak 
\fi 
\H@refstepcounter{Hfootnote}%
\hyper@makecurrent{Hfootnote}%
\hyper@linkstart{link}{\@currentHref}%
\@makefnmark 
\hyper@linkend 
\ifFN@pp@towrite 
\FN@pp@writetemp 
\FN@pp@towritefalse 
\fi 
\FN@mf@prepare 
\ifhmode\spacefactor\@x@sf\fi 
\relax%

\makeatother 

\theoremstyle{plain}
\newtheorem{theorem}{Theorem}[section]

\newaliascnt{proposition}{theorem}
\newtheorem{proposition}[proposition]{Proposition}
\aliascntresetthe{proposition}

\newaliascnt{lemma}{theorem}
\newtheorem{lemma}[lemma]{Lemma}
\aliascntresetthe{lemma}

\newaliascnt{corollary}{theorem}
\newtheorem{corollary}[corollary]{Corollary}
\aliascntresetthe{corollary}

\newaliascnt{conjecture}{theorem}

\aliascntresetthe{conjecture}

\newaliascnt{property}{theorem}

\aliascntresetthe{property}

\newaliascnt{claim}{theorem}
\newtheorem{claim}[claim]{Claim}
\aliascntresetthe{claim}

\theoremstyle{definition}
\newaliascnt{definition}{theorem}
\newtheorem{definition}[definition]{Definition}
\aliascntresetthe{definition}

\newaliascnt{example}{theorem}

\aliascntresetthe{example}

\newaliascnt{question}{theorem}

\aliascntresetthe{question}

\newaliascnt{notation}{theorem}

\aliascntresetthe{notation}

\theoremstyle{remark}
\newaliascnt{remark}{theorem}
\newtheorem{remark}[remark]{Remark}
\aliascntresetthe{remark}

\newaliascnt{observation}{theorem}

\aliascntresetthe{observation}

\makeatletter
\newtheorem*{rep@theorem}{\rep@title}
\newcommand{\newreptheorem}[2]{%
\newenvironment{rep#1}[1]{%
 \def\rep@title{#2 \ref{##1}}%
 \begin{rep@theorem}}%
 {\end{rep@theorem}}}
\makeatother

\newcommand{\eps}{{\varepsilon}}
\newcommand{\sys}{{\rm sys}}
\newcommand{\ii}{{\mathit i}}
\newcommand{\Sys}{{\mathfrak S}}
\newcommand{\R}{{\mathbb R}}
\newcommand{\D}{{\mathbb D}}
\newcommand{\Hyp}{{\mathbb H}}
\newcommand{\N}{{\mathbb N}}
\newcommand{\F}{{\mathbb F}}
\newcommand{\id}{{\mathrm{id}}}
\newcommand{\TT}{{\mathbb T}}
\newcommand{\K}{{\mathbb K}}
\newcommand{\Z}{{\mathbb Z}}
\newcommand{\C}{{\mathbb C}}
\newcommand{\GG}{{\mathbb G}}
\newcommand{\GGS}{{\mathbb G}(S)}
\newcommand{\CC}{{\mathcal C}}
\newcommand{\BB}{{\mathcal B}}
\newcommand{\OO}{{\mathcal O}}
\newcommand{\arc}{{\mathcal A}}
\newcommand{\DD}{{\mathcal D}}
\newcommand{\Sp}{{\mathbb S}}
\newcommand{\G}{{\mathcal G}}
\newcommand{\T}{{\mathcal T}}
\newcommand{\M}{{\mathcal M}}
\newcommand{\PSL}{{\rm PSL}}
\newcommand{\area}{{\rm area}}
\newcommand{\arcsinh}{{\,\rm arcsinh}}
\newcommand{\arccosh}{{\,\rm arccosh}}
\newcommand{\diam}{{\rm diam}}
\newcommand{\dias}{{\rm dias}}
\newcommand{\numS}{{\rm Kiss}(S)}
\newcommand{\length}{{\rm length}}
\newcommand{\measure}{\lambda}
\newcommand{\bigmeasure}{\lambda_{\vec{k}}}
\newcommand{\1}{\mathbf{1}}
\newcommand{\kk}{\Bbbk}
\newcommand{\PP}{\mathcal P}
\newcommand{\FF}{\mathcal F}
\newcommand{\II}{\mathcal I}
\newcommand{\cf}{{\it cf.}}
\newcommand{\ie}{{\it i.e.}}
\newcommand{\curves}{\mathcal H}

\newcommand{\Li}{{\rm Li}}
\newcommand{\Ree}{\operatorname{Re}}
\newcommand{\mcg}{{\rm Mod}}
\newcommand{\Log}{{\rm Log}}
\newcommand{\Arg}{{\rm Arg}}

\def\bull{\vrule height 1ex width .8ex depth -.2ex}
\renewcommand{\labelitemi}{\bull}

\renewcommand{\footnoterule}{\vfill \rule{0.25\textwidth}{0.5pt} \vspace{4pt}}

\linespread{1.25}

\sectionfont{\large \bfseries}
\subsectionfont{\normalsize}

\setlength{\parindent}{0pt}
\setlength{\parskip}{6pt}

\long\def\symbolfootnote[#1]#2{\begingroup%
\def\thefootnote{\fnsymbol{footnote}}\footnote[#1]{#2}\endgroup}

\def\blfootnote{\xdef\@thefnmark{}\@footnotetext}

\begin{document}

{\Large \bfseries Graded Bridgeman dilogarithm identities on hyperbolic surfaces}

{\large Ara Basmajian, 
Nhat Minh Doan,
Hugo Parlier,
 and Ser Peow Tan}


\begingroup
 \renewcommand\thefootnote{}%
\endgroup

{\bf Abstract.} We establish graded versions of Bridgeman’s dilogarithm identity for hyperbolic cone surfaces, including surfaces with only cusps and cone points, and provide applications to the study of orthogeodesics.

\begin{minipage}{\textwidth}
\begingroup
\hypersetup{linkcolor=black}
\tableofcontents
\endgroup
\end{minipage}


\section{Introduction}\label{sec:intro}

In the study of hyperbolic surfaces, a class of universal length identities has emerged: infinite sums of explicit functions of lengths of geodesics that remain constant, either as a topological invariant or as the geodesic/horocyclic boundary length, across every hyperbolic metric in a given moduli space. Well-known examples include McShane’s identity for cusped surfaces \cite{mcshane1998simple}, extended by Mirzakhani to bordered surfaces (hyperbolic surfaces with geodesic boundary) \cite{mirzakhani2007simple} and by Tan–Wong–Zhang to surfaces with boundary and/or cone points \cite{tan2006generalizations}, Basmajian and Bridgeman’s orthospectrum identities \cite{basmajian1993orthogonal,bridgeman2011orthospectra} for surfaces with geodesic boundary, and Luo-Tan's identity \cite{luo2014dilogarithm} for surfaces with or without geodesic boundary. All of these identities arise from distinct geometric decompositions that reflect different behaviors of the geodesic flow on hyperbolic surfaces (see \cite{bridgeman2016identities}). In a slightly different direction, Wang and Xue \cite{WangXue-CPAM-2025} recently discovered a new type of identity that relates the length spectrum of a surface to that of a closely related surface obtained by adding cusps to the original surface via uniformization. 

Basmajian’s identity states that, for a hyperbolic surface $X$ with (at least one) geodesic boundary,
$$
\sum_{\mu}2\ln\coth\Bigl(\tfrac{\ell(\mu)}2\Bigr)
=\ell(\partial X),
$$
where the sum runs over all orthogeodesics $\mu$ on $X$ and $\ell(\mu)$ is its length. Crucially, Basmajian’s identity requires at least one closed geodesic boundary. On a surface with only cusps it becomes trivial, since every orthogeodesic then has infinite length, and unlike for McShane's identity, replacing geodesic length with horocyclic length produces strict inequalities.

Basmajian, Parlier and Tan \cite{basmajian2025prime} overcame this obstacle in generalizing Basmajian’s identity by developing a unified framework that treats geodesic boundaries, cusps and cone points all at once. One of the key ideas is putting grades on each boundary element (cusp, cone, or geodesic boundary component), and considering natural collars around each boundary element associated to the grade. Removing these collar neighborhoods yields a compact subsurface, the so-called \emph{concave core}. The boundary of the concave core can then be partitioned into disjoint sets indexed by \emph{prime orthogeodesics}, that is orthogeodesic arcs contained entirely within the concave core that satisfy a specific primality condition, together with a complementary set of zero measure. Their family of identities varies continuously with the type of boundary: as one transitions from a hypercycle (an equidistant curve to a geodesic boundary) to a horocycle around a cusp and further to a circle around a cone point, the identity interpolates between Basmajian’s identity and its new cusped/conical variants. 

From a different viewpoint, Bridgeman \cite{bridgeman2011orthospectra} proved a striking dilogarithmic identity linking the orthogeodesic spectrum of a hyperbolic surface with boundary to the volume of its unit tangent bundle. Concretely, if $X$ has at least one closed geodesic boundary, then
$$
\sum_{\mu}\mathcal{L}\left(\frac{1}{\cosh^2(\tfrac{\ell(\mu)}{2})}\right)
=\frac{\pi^2}{2}\left|\chi(X)\right|,
$$
where the sum runs over all orthogeodesics $\mu$ of length $\ell(\mu)$, and $\mathcal{L}(.)$ is the Rogers dilogarithm \cite{rogers1907function} defined on the principal branch \(\mathbb C\setminus ((-\infty,0]\cup[1,+\infty))\) by
$$
\mathcal{L}(z):=\operatorname{Li}_2(z)+\tfrac12\Log(z)\,\Log(1-z)=-\frac12\int_{0}^{z}\left(\frac{\Log(1-t)}{t}+\frac{\Log(t)}{1-t}\right)\,dt.
$$
 By the ergodicity of the geodesic flow, the geodesic‐boundary assumption ensures that almost every unit tangent vector, when followed forward and backward, traces out a geodesic segment that meets the boundary in finite time, and which is homotopic relative to the boundary to a unique orthogeodesic $\mu$. This allows one to partition the unit tangent bundle $T^1(X)$ into regions labeled by each $\mu$, each region having Liouville volume given by a universal dilogarithm function of $\ell(\mu)$. Since Gauss–Bonnet gives $\mathrm{Vol}(T^1(X))=4\pi^2|\chi(X)|$, summing these contributions recovers the stated identity.

 Luo and Tan \cite{luo2014dilogarithm} extended Bridgeman's idea to every hyperbolic surface with or without geodesic boundary by stopping each geodesic as soon as it either self-intersects (forming a loop) or first hits the boundary. Except on a measure-zero set, applying this rule to both the forward and backward rays from any unit tangent vector produces a “spine” graph, which sits inside a unique embedded pair of pants or one-holed torus. Finally, by grouping all vectors whose spines lie in the same subsurface and summing their Liouville volumes over every such subsurface, one obtains the Luo-Tan identity. 

In this paper, we adapt the ideas of \cite{basmajian2025prime} to prove graded versions of the Bridgeman identity, which hold even when the surface has only cusps and/or cone points. As in \cite{basmajian2025prime}, for a hyperbolic cone surface $X$, we put grades $k_i\in \mathbb N \cup \{\infty\}$ on each boundary element. Removing collars adapted to the grades around each of the boundary elements from $X$ gives the concave core $V_{\vec k}(X)$. Associated with each {\it prime orthogeodesic} $\mu \subset V_{\vec k}(X)$ is a {\it maximal immersed pair of pants} $P_{\vec k}(\mu)$. We show how a generic unit vector $v \in T^1(X)$ is associated to exactly one of these pairs of pants by showing how $v$ generates a geodesic arc which is immersed in a unique maximal immersed pair of pants. This arc is obtained by exponentiating $v$ in both directions until just before it reaches the deepest point of a collar, or hits a boundary geodesic. This gives a decomposition of a full measure subset of $T^1(X)$ into disjoint subsets indexed by the prime orthogeodesics. The measure of each of these subsets can then be computed using ideas from \cite{bridgeman2011orthospectra} and \cite{luo2014dilogarithm}. Note that all the vectors tangent to the geodesic arc generated by $ v$ generate the same geodesic arc, and hence lie in the same subset of the decomposition, and that the measures of these sets depend on the geometry of the immersed pair of pants $P_{\vec k}(\mu)$, and not just on the length of the prime orthogeodesic $\ell(\mu)$.

In extending the identity to hyperbolic cone surfaces, we address two main difficulties. 
First, we obtain an explicit closed form (via Rogers dilogarithms) for the Lasso function associated to a cone point by evaluating the double integral
$$I(a,b) = \int_{a}^{b}\int_{-1/y}^{-b} \frac{\ln\left|\frac{(x-a)(x-b)(1+y^2)}{(y-a)(y-b)(1+x^2)}\right|}{(y-x)^2}\,dx\,dy,$$
Second, we prove that the set of generic vectors has full Liouville measure even when non-orbifold cone points are present. For context, while the ergodicity of the geodesic flow is classical for finite-area hyperbolic $2$-orbifolds \cite{hopf1939statistik,ansov1969geodesic,hopf1971ergodic,Nicholls1989,katok1995introduction}, no general ergodicity theorem is known for hyperbolic cone surfaces with non-orbifold cone points. The geometry of hyperbolic cone surfaces has been studied extensively, from the classical theory of conical metrics \cite{mcowen1988point,Troyanov1991} to more recent work on their moduli spaces and Weil–Petersson volumes \cite{do2009weil,anagnostou2022volumes}, as well as related developments in mathematical physics \cite{bonsante2009ads,witten2020matrix,almheiri2020replica}.

\subsection{Main results}
Our main result is the identity stated in \cref{thm:AbstractIdentity} which gives a decomposition of a full measure subset of $T^1(X)$ into subsets indexed by $\vec k$-prime orthogeodesics and the formulas for the measures of each of these subsets in terms of the geometry of the immersed pair of pants $P_{\vec k}(\mu)$. As a first application, we prove a counting result for orthogeodesics (see \cref{thm:Counting}), which behaves roughly like $e^L/L$ as in Huber's theorem (see \cite{buser2010geometry}). We then deduce a non-domination theorem for the spectrum of $\vec k$-prime orthogeodesics on hyperbolic surfaces with only cusps (see \cref{Nondomination}), which is an analog of Thurston’s theorem \cite{thurston1998minimal} on the marked simple length spectrum of closed hyperbolic surfaces.

Before stating the main theorem, we introduce the setting, notation, and definitions, partly following \cite{basmajian2025prime}. 
Let $X$ be a hyperbolic cone surface obtained from a complete hyperbolic cone surface with exactly $n>0$ boundary elements (simple closed geodesic boundary components, cone points, or cusps) by deleting all cone points. 
We continue to refer to these deleted cone points as the boundary elements of $X$, and we orient each of them so that $X$ lies to their right.
For a boundary element $\delta$ of $X$, define its boundary parameter (length for geodesic boundaries, angle for cone points) by
$$
b(\delta)=
\begin{cases}
\ell(\delta), & \text{if $\delta$ is a geodesic boundary component},\\
\theta(\delta), & \text{if $\delta$ is a cone point}.
\end{cases}
$$
Cusps are regarded either as cone points of angle $0$ or as geodesic boundary components of length $0$.
Assign a grade $k_i\in\overline{\mathbb{N}}=\{1,2,3,\ldots\}\cup\{\infty\}$ to each boundary element $i \in \{1,2,\ldots,n\}$. A choice of grades is a \emph{grading} denoted
\begin{equation}\label{eq:Grading}
 \vec{k}=(k_1,k_2,\ldots,k_n).
\end{equation}
\begin{definition}[Admissible pair]\label{dfn:AdmissiblePair}
The pair $(X,\vec k)$ is \emph{admissible} if:
\begin{enumerate}
 \item For every index $i$ such that the $i$-th boundary element is a cone point, the corresponding angle $\theta_i$ satisfies $0\le \theta_i \le \pi/{k_i}$.
 \item There exists an index $i$ such that either the $i$th boundary element is not a cone point, or it is a cone point with angle $\theta_i$ satisfies $0\le \theta_i<\pi/k_i$.
 \item $X$ is not a pair of pants whose grading $\vec k$ is a permutation of $(1,1,m)$ for any $m\in\overline{\mathbb{N}}$.
\end{enumerate}
\end{definition}

Condition (1) guaranties well-defined collars around each boundary element (see \cref{dfn:kNaturalCollar}); (2) excludes the degenerate case with no collars and no geodesic boundary components; and (3) removes the trivial one-term case of \cref{thm:AbstractIdentity} and ensures the existence of a model surface (see \cref{dfn:modelSurface}). In this paper, we consider only admissible pairs $(X,\vec k)$.

The \emph{$k_i$-th} natural collar of the $i$th boundary element is denoted by $C_{k_i}$. Roughly speaking, it is the smallest collar about the $i$th boundary element such that any geodesic arc forming a $k_i$ loop around it intersects the collar (see \cref{dfn:kNaturalCollar}). The \emph{$\vec{k}$-concave core} of $X$ is defined as
\begin{equation}\label{eq:Concave}
V_{\vec{k}}(X) = X \setminus \bigcup_{i=1}^{n}C_{k_i}.
\end{equation}
More generally, removing disjoint collars of any non-negative radii yields a concave core of $X$. A (truncated) orthogeodesic in the $\vec{k}$‑concave core $V_{\vec{k}}(X)$ is an immersed, oriented geodesic segment $\eta \subset V_{\vec{k}}(X)$
that meets the boundary $\partial V_{\vec{k}}(X)$ orthogonally at both endpoints. Each such orthogeodesic extends uniquely to an oriented \emph{$\vec{k}$‑orthogeodesic} on $X$. Let $\OO_{\vec{k}}(X)$ and $\overline{\OO}_{\vec{k}}(X)$ be the sets of all oriented and unoriented $\vec{k}$-orthogeodesics on $X$, respectively.

\begin{definition}[Immersed pair of pants $P_{\vec{k}}(\mu)$]\label{dfn:immersedPants}
 Let $\mu\in\mathcal{O}_{\vec{k}}(X)$ be a $\vec{k}$-orthogeodesic oriented from the boundary element $\alpha$ (grade $k_\alpha<\infty$) to $\beta$ (grade $k_\beta<\infty$). 
We assume $k_\alpha,k_\beta<\infty$; the degenerate cases are treated in \cref{prop:functionf}. 
Under \cref{dfn:AdmissiblePair}(3), the loop
\begin{equation}\label{eq:Gamma}
\alpha^{k_\alpha} * \mu * \beta^{k_\beta} * \mu^{-1}
\end{equation}
represents a non-peripheral free homotopy class on $X$. When this class is realized by a unique closed geodesic on $X$, we denote this geodesic by $\gamma_{\vec{k}}(\mu)$. 
 The triple $(\alpha^{k_\alpha}, \beta^{k_\beta}, \gamma_{\vec{k}}(\mu))$ is then the boundary of a \emph{$\vec k$-immersed pair of pants} $P_{\vec{k}}(\mu)$. Equivalently, we may specify a \emph{pre-immersed pair of pants} $\widetilde{P}_{\vec{k}}(\mu)$ and a local isometry
\begin{equation}\label{eq:PrePants}
 \varphi\colon \widetilde{P}_{\vec{k}}(\mu) \longrightarrow X 
\end{equation}
such that $\varphi(\widetilde{P}_{\vec{k}}(\mu)) = P_{\vec{k}}(\mu)$. The surface $\widetilde{P}_{\vec{k}}(\mu)$ is a hyperbolic pair of pants whose boundary components, denoted
\[
\widetilde{\alpha}^{k_{\alpha}},\qquad
\widetilde{\beta}^{k_{\beta}},\qquad
\widetilde{\gamma}_{\vec{k}}(\mu),
\]
map respectively to $\alpha^{k_\alpha}$, $\beta^{k_\beta}$, and $\gamma_{\vec{k}}(\mu)$ under $\varphi$ and have boundary (length or angle) parameters
\[
k_\alpha\,b(\alpha),\qquad k_\beta\,b(\beta),\qquad \ell\big(\gamma_{\vec{k}}(\mu)\big).
\]
See \cref{fig:merged} for the case $k_{\alpha}=k_{\beta}=2$ with geodesic boundary components $\alpha$ and $\beta$.

\end{definition}
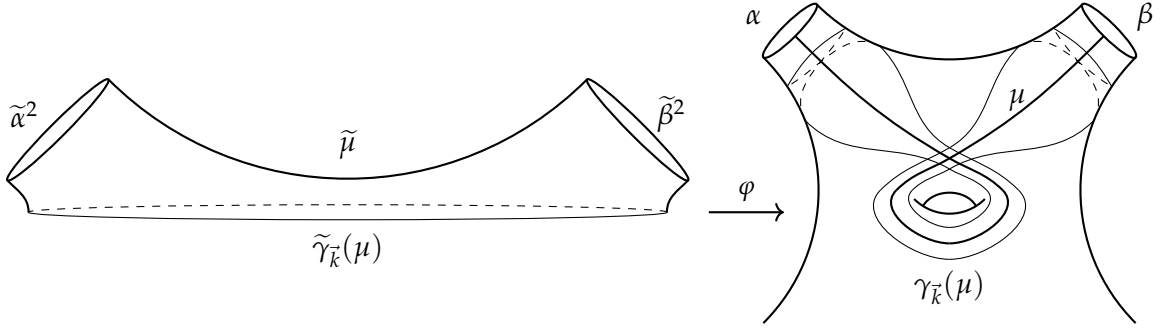
\begin{figure}[h!]
 \centering
 \begin{tikzpicture}
 
\begin{scope}[yshift=-1cm]
 

\node at (0,0.9) {$\widetilde{\mu}$};

\node at (4.3,1.3) {$\widetilde{\beta}^{2}$};
 
\node at (-4.3,1.3){$\widetilde{\alpha}^{2}$};
 
\node at (0,-0.5) {$\widetilde{\gamma}_{\vec{k}}(\mu)$};

\begin{scope}[rotate=45, scale=0.5]
 \draw[thick] (-3.9,7) ellipse (1.9 and 0.2);
 \draw[thick] (7,-3.9) ellipse (0.2 and 1.9);

 \draw[thin,dashed,rotate around={45:(0,0)}]
 (0,-8.5) 
 arc[x radius=0.2,y radius=8.5, start angle=-90, end angle=90];

 \draw[thin,rotate around={45:(0,0)}]
 (0,8.5) 
 arc[x radius=0.2,y radius=8.5, start angle=90, end angle=270];


 \coordinate (D1) at (-2,7);
 \coordinate (D2) at (7,-2);
 \coordinate (D3) at (-5.8,7);
 \coordinate (D4) at (-6,6);
 \coordinate (D5) at (6,-6);
 \coordinate (D6) at (7,-5.8);


 \draw[thick] (D1) to[out=-90,in=180] (D2);
 
 \draw[thick] (D3) to[out=-90,in=45] (D4);
 \draw[thick] (D5) to[out=45,in=-180] (D6);
\end{scope}
\end{scope}

 \begin{scope}[xshift=5.8cm,yshift=-1cm]
 \draw[->, thick] (-1,0) -- (0,0) node[midway, above]{\small $\varphi$};
 \end{scope}

\begin{scope}[xshift=8cm, scale=1]

\node at (0.9,0.5) {$\mu$};

\node at (2.6,1.6) {$\beta$};
 
\node at (-2.6,1.6){$\alpha$};
 
\node at (0,-2) {$\gamma_{\vec{k}}(\mu)$};

\begin{scope}[rotate=45, scale=0.5]
 
 \coordinate (A3) at (-2,0);
 \coordinate (A4) at (-2,-2);
 \coordinate (A5) at (0,-2);
 \coordinate (A6) at (-1,4.8);
 \coordinate (A7) at (4.8,-1);
 \coordinate (A2) at (-0.4,-0.4);
 
\foreach \i in {}{
 \node[font=\tiny, inner sep=1pt] at (A\i) {$A_{\i}$};
}

\draw[thick] (A3) to[out=-180,in=135] (A4);
 
\draw[thick] (A4) to[out=-45,in=-90] (A5);

\draw[thick] (A2) to[out=105,in=-90] (A6);

\draw[thick] (A2) to[out=-15,in=180] (A7);
 
\draw[thick] (A3) to[out=0,in=165] (A2);

\draw[thick] (A5) to[out=90,in=-75] (A2);

 
 \coordinate (B4) at (-2.3,-2.3);
 \coordinate (C4) at (-1.7,-1.7);
 \coordinate (B3) at (-2.3,0.3);
 \coordinate (C3) at (-1.7,-0.3);
 \coordinate (B5) at (0.3,-2.3);
 \coordinate (C5) at (-0.3,-1.7);
 \coordinate (B2) at (-0.2,-0.2);
 \coordinate (C2) at (-0.6,-0.6);
 \coordinate (C1) at (-2.4,3);
 \coordinate (B1) at (0.4,3);
 \coordinate (C0) at (3,-2.4);
 \coordinate (B0) at (3,0.4);
 \coordinate (C6) at (0.1,4);
 \coordinate (B6) at (-2.1,4);
 \coordinate (C7) at (4,0.1);
 \coordinate (B7) at (4,-2.1);
 
 \foreach \i in {}{
 \node[font=\tiny, inner sep=1pt] at (B\i) {$B_{\i}$};
} 

 \foreach \i in {}{
 \node[font=\tiny, inner sep=1pt] at (C\i) {$C_{\i}$};
} 

\draw[thin] (B3) to[out=-180,in=135] (B4);

\draw[thin] (C3) to[out=-180,in=135] (C4);
 
\draw[thin] (B4) to[out=-45,in=-90] (B5);

\draw[thin] (C4) to[out=-45,in=-90] (C5);

\draw[thin] (B3) to[out=0,in=165] (B2);

\draw[thin] (C3) to[out=0,in=165] (C2);

\draw[thin] (B5) to[out=90,in=-75] (B2);

\draw[thin] (C5) to[out=90,in=-75] (C2);

\draw[thin] (C2) to[out=105,in=-90] (C1);

\draw[thin] (B2) to[out=105,in=-90] (B1);

\draw[thin] (B2) to[out=-15,in=180] (B0);

\draw[thin] (C2) to[out=-15,in=180] (C0);

\draw[dashed] (B1) to[out=105,in=-15] (B6);

\draw[dashed] (C1) to[out=75,in=-165] (C6);

\draw[thin] (B6) to[out=15,in=165] (C6);

\draw[dashed] (B0) to[out=-15,in=105] (B7);

\draw[dashed] (C0) to[out=15,in=-105] (C7);

\draw[thin] (B7) to[out=75,in=-75] (C7);
 
 \draw[thick] (-1.83,-0.5) to[out=-90,in=180] (-0.5,-1.83);
 \draw[thick] (-1.8,-0.8) to[out=15,in=75] (-0.8,-1.8);
 
 \draw[thick] (0,5) to[out=-90,in=180] (5,0);
 \draw[thick] (-2,5) to[out=-90,in=0] (-7,0);
 \draw[thick] (0,-7) to[out=90,in=-180] (5,-2);
 
 \draw[thick] (-1,5) ellipse (1 and 0.2);
 \draw[thick] (5,-1) ellipse (0.2 and 1);

 \end{scope}
 \end{scope}

 \end{tikzpicture}
 \caption{Illustration of the isometric immersion $\varphi:\widetilde{P}_{\vec{k}}(\mu)\rightarrow X$ when $k_{\alpha} = k_{\beta} = 2$.}
 \label{fig:merged}
\end{figure}

Note that reversing the orientation of \(\mu\) yields the same immersed pair of pants. Both $\mu$ and $-\mu$ are called the {\it spine} of $P_{\vec{k}}(\mu)$. 

\begin{definition}[$\vec{k}$-maximal immersed pairs of pants and $\vec{k}$-prime orthogeodesics]\label{dfn:kPrimeOrthogeodesic}
Let $\mu'\in\overline{\mathcal{O}}_{\vec{k}}(X)$, and suppose the corresponding immersed pair of pants $P_{\vec{k}}(\mu')$ in $X$ exists. 
We say that $P_{\vec{k}}(\mu')$ is \emph{$\vec{k}$-maximal} if there is no other orthogeodesic $\mu \in \overline{\mathcal{O}}_{\vec{k}}(X)$ with $\mu\neq\mu'$ such that $P_{\vec{k}}(\mu')$ isometrically immerses into $P_{\vec{k}}(\mu)$. 
In that case, we call $\mu'$ a \emph{$\vec{k}$-prime orthogeodesic}.

Let $\mathcal{O}'_{\vec{k}}(X)$ and $\overline{\mathcal{O}}'_{\vec{k}}(X)$ denote the sets of $\vec{k}$-prime oriented and unoriented orthogeodesics, respectively. 
Under \cref{dfn:AdmissiblePair}(3), the map $\mu'\mapsto P_{\vec{k}}(\mu')$ is a bijection from $\overline{\mathcal{O}}'_{\vec{k}}(X)$ onto the set of $\vec{k}$-maximal immersed pairs of pants in $X$.

\end{definition}

\begin{remark}\label{rem:isometricImmersion}

 The phrase "$P_{\vec{k}}(\mu')$ is isometrically immersed in $P_{\vec{k}}(\mu)$" means that there exists an isometric immersion $\psi$
making the following diagram commute:

\begin{center}
\begin{tikzcd}[column sep=large, row sep=large]
\widetilde{P}_{\vec{k}}(\mu') 
 \arrow[rr, "\psi"] 
 \arrow[dr,swap,"\varphi'"] 
 && \widetilde{P}_{\vec{k}}(\mu) 
 \arrow[dl, "\varphi"] \\
& X &
\end{tikzcd}
\end{center}

Here, $\varphi$ and $\varphi'$ are the isometric immersions into $X$ corresponding respectively to $P_{\vec{k}}(\mu)$ and $P_{\vec{k}}(\mu')$. We note that $P_{\vec{k}}(\mu')\subset P_{\vec{k}}(\mu)$ as sets if $P_{\vec{k}}(\mu')$ is isometrically immersed in $P_{\vec{k}}(\mu)$. Conversely, if $P_{\vec{k}}(\mu')\subset P_{\vec{k}}(\mu)$ as sets and $\varphi'_{*}(\pi_{1}(\widetilde P_{\vec k}(\mu'),\widetilde{x}'_0))\subset \varphi_{*}(\pi_{1}(\widetilde P_{\vec k}(\mu),\widetilde{x}_0))$ as subgroups of $\pi_1(X,x_0)$, where $x_0 = \varphi(\widetilde{x}_0) = \varphi'(\widetilde{x}'_0)$, then $P_{\vec{k}}(\mu')$ is isometrically immersed in $P_{\vec{k}}(\mu)$. Under this isometric immersion relation, the collection of all immersed pairs of pants forms a poset (partially ordered set) of infinite height. See \cite{basmajian2025prime} for other equivalent definitions of a $\vec{k}$-prime orthogeodesic; in particular, see \cref{prop:prime-realization} for its realization in the model surface.
\end{remark}

\begin{definition}[Non-generic vectors]\label{dfn:NongenericVector}
 Let $T_{1}(X)$ be the unit tangent bundle of $X$. For each 
$v\in T_{1}(X)$ we follow the geodesic determined by $v$ in both forward and backward direction at equal speed. A vector $v \in T^1(X)$ is called \emph{non-generic} if $v \in B_{1}\cup B_{2}\cup B_{3}$, where
 \begin{enumerate}
 \item $B_1$ is the set of vectors whose base points lie in a collar and are tangent to a horocycle, hypercycle, or circle associated to the corresponding boundary element (horocycle for a cusp, hypercycle for a geodesic boundary, circle for a cone point);
 
 \item $B_2$ is the set of those vectors whose base points lie in a collar, and for which the direction which first exits the collar stays in the concave core $V_{\vec{k}}(X)$ forever;
 
 \item $B_3$ is the set of those vectors whose base points lie in the concave core, and for which at least one direction stays in the concave core $V_{\vec{k}}(X)$ forever.
 \end{enumerate}
\end{definition}
The set $B_1$ has Liouville measure zero. When $X$ is a hyperbolic orbifold, ergodicity of the geodesic flow on $T^1X$ implies that $B_2$ and $B_3$ also have Liouville measure zero. More generally, we will show in \cref{sec:Nongeneric} (\cref{pro:WeakErgodicity}) that the same conclusion holds for hyperbolic cone surfaces.

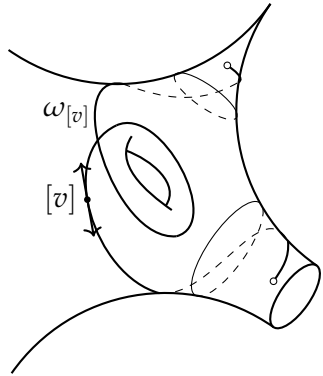
\begin{figure}[h!]
 \centering
\begin{tikzpicture}[scale=0.5, rotate=-35]

 \coordinate (A1) at (0,3);
 \coordinate (A2) at (0.5,2.8);
 \coordinate (A3) at (-2,1);
 \coordinate (A4) at (-1,-2.06);
 \coordinate (A5) at (2,-3);
 \coordinate (A6) at (4,0.1);
 \coordinate (A7) at (4.3,-1);
 \node[left] at (A4) {$[v]$};

 \node[left] at (-2,0) {$\omega_{[v]}$};

 \draw[thick] (0,5) to[out=-90,in=180] (5,0);
 \draw[thick] (0,5) to[out=-90,in=0] (-5,0);
 \draw[thick] (0,-7) to[out=90,in=-180] (5,-2);

 \draw[thick] (A1) to[out=0, in=65] (A2);
 \draw[fill=white,draw=black] (A1) circle (2.5pt);

 \draw[dashed] (A3) to[out=15, in=-115] (A2);
 \draw[thick] (A3) to[out=-155,in=220] (1.5,-1.5);
 \draw[thick] (1.5,-1.5) to[out=50,in=60] (-1.5,0);
 \draw[thick] (-1.5,0) to[out=-120,in=200] (A5);
 \draw[dashed] (A5) to[out=30, in=155] (A6);
 \draw[thick] (A6) to[out=-50, in=90] (A7);
 \draw[fill=white,draw=black] (A7) circle (2.5pt);

 \draw[thick] (-1,0) to[out=-90,in=180] (1,-1);
\draw[thick] (-0.9,-0.4) to[out=15,in=90] (0.8,-1);

 \draw[->,thick] (A4) -- ++(135:1cm);
 \draw[->,thick] (A4) -- ++(-45:1cm);
 \fill (A4) circle (2.5pt);

 \draw[thick] (5,-1) ellipse (0.4 and 1);
 \draw[thin] (2.8,-1) ++(90:1.5) arc(90:270:0.5 and 1.5);
 \draw[thin,dashed] (2.8,-1) ++(-90:1.5) arc(-90:90:0.5 and 1.5);

 \draw[thin,dashed] (0,2) ++(180:1 and 0.3) arc(180:360:1 and 0.3);
 \draw[thin] (0,2) ++(0:1 and 0.3) arc(0:180:1 and 0.3);

\end{tikzpicture}
 \caption{Illustration of the geodesic path \(\omega_{[v]}\) determined by a generic vector \(v\in T^1(V_{\vec{k}}(X))\).}

 \label{fig:genericVector}
\end{figure}

\begin{definition}[Generic vector and its associated geodesic path]\label{dfn:genericVector}
A vector $v\in T^1(X)\setminus (B_{1}\cup B_{2}\cup B_{3})$ is called \emph{generic}. Writing $[v]=\{\pm v\}$ for a generic vector $v$, its associated immersed geodesic path $\omega_{[v]}$ is constructed as follows:

\begin{enumerate}
 \item If $\operatorname{base}(v)\in V_{\vec k}(X)$, then exponentiate $v$ in the forward direction until it reaches the boundary of $V_{\vec{k}}(X)$, and continue the geodesic into the boundary collar until just before it hits its deepest point in the collar. Do the same for the backward direction (see \cref{fig:genericVector}). This defines $\omega_{[v]}$ in this case.
 \item If $\operatorname{base}(v)$ lies in a boundary collar, then in one direction the geodesic reaches the deepest point inside that collar, while in the opposite direction it exits the collar (i.e., enters the concave core) and continues until just before it reaches the deepest point in another boundary collar (which might be the same as the original). This defines $\omega_{[v]}$ in this case.
\end{enumerate}
\end{definition}

The following key technical proposition establishes a correspondence between generic vectors and maximal immersed pairs of pants, which underlies our identity. Recall that for an orthogeodesic $\mu \subset V_{\vec k}(X)$, we can define a $\vec k$-immersed pair of pants $P_{\vec k}(\mu)$ which has $\mu$ as its spine (see \cref{dfn:immersedPants}). There is a partial ordering on the set of $\vec k$-immersed pairs of pants defined by isometric immersion, and $\mu$ is said to be prime if $P_{\vec k}(\mu)$ is maximal with respect to this partial ordering (see \cref{dfn:kPrimeOrthogeodesic}).

\begin{proposition}\label{prop:uniqueness}
For each generic vector \(v \in T^1(X)\), there exists a unique $\vec k$-immersed pair of pants 
$P_{\vec{k}}(\mu),$ where \(\mu\) is a \(\vec k\)-prime orthogeodesic, such that 
$ \omega_{[v]}$ is isometrically immersed in $P_{\vec{k}}(\mu)$.
\end{proposition}
We will prove \cref{prop:uniqueness} in \cref{sec:abstract} by leveraging the model surface and its punctured orbifold cover (see \cref{dfn:modelSurface} and \cref{dfn:punctureOrbifoldCover}), as introduced in \cite{basmajian2025prime}.

\begin{remark}\label{rem:isometricImmersionArc}
 Similar to \cref{rem:isometricImmersion}, the phrase "$\omega_{[v]}$ is isometrically immersed in $P_{\vec{k}}(\mu)$" in \cref{prop:uniqueness} means that there exists an isometric immersion $\psi_I$ making the following diagram commute:
\begin{center}
\begin{tikzcd}[column sep=large, row sep=large]
I 
 \arrow[rr, "\psi_I"] 
 \arrow[dr, "\varphi_I"'] 
 && \widetilde{P}_{\vec{k}}(\mu) 
 \arrow[dl, "\varphi"] \\
& X &
\end{tikzcd}
\end{center}
Here, $I$ is an open geodesic segment in $\mathbb{H}^2$ of length $\ell(\omega_{[v]})$, and $\varphi_I$ and $\varphi$ denote the isometric immersions into $X$ corresponding to $\omega_{[v]}$ and $P_{\vec{k}}(\mu)$, respectively.

\end{remark}

 If $v$ is a generic vector, then every unit tangent vector $u$ tangent to $\omega_{[v]}$ is also generic since $\omega_{[u]}=\omega_{[v]}$. The non-generic set $B_1\cup B_2\cup B_3$ (see \cref{dfn:NongenericVector}) has Liouville measure zero (see \cref{pro:WeakErgodicity}); combined with \cref{prop:uniqueness}, this yields the following abstract identity:

\begin{theorem}\label{thm:AbstractIdentity} 
Let $(X,\vec k)$ be an admissible pair (see \cref{dfn:AdmissiblePair}), and let $\overline{\OO}'_{\vec k}(X)$ denote the set of $\vec k$-prime unoriented orthogeodesics (see \cref{dfn:kPrimeOrthogeodesic}). Then
\[
\sum_{\mu\in \overline{\OO}'_{\vec{k}}(X)} h\bigl(P_{\vec{k}}(\mu)\bigr)
= \operatorname{Vol}\bigl(T^1(X)\bigr),
\]
where $h\bigl(P_{\vec{k}}(\mu)\bigr)$ is the Liouville volume of the set of generic $v\in T^1(X)$ such that $\omega_{[v]}$ is isometrically immersed in $P_{\vec{k}}(\mu)$. 
Moreover, if for every index $i$ such that the $i$-th boundary element of $X$ is a cone point, the corresponding angle $\theta_i$ satisfies $ \theta_i \le \pi/2k_i$, then $h\bigl(P_{\vec{k}}(\mu)\bigr)$ is given explicitly in \cref{prop:functionf}.

\end{theorem}
\begin{remark}\label{rem:SpecialFunctionf}

\begin{enumerate}

 \item In \cref{thm:AbstractIdentity},
 the cone–angle bound prevents overlap in the geometric decomposition used to compute $h$; see \cref{lem:OutQ,subsec:SmallConeAngle}.
 
 \item The proof of our identity parallels that of Luo–Tan \cite{luo2014dilogarithm}, replacing embedded pairs of pants with immersed ones. By \cref{prop:uniqueness}, each generic vector lies in a unique $\vec k$-maximal immersed pair of pants $P_{\vec k}(\mu)$, associated to a unique $\vec k$-prime orthogeodesic $\mu\in\overline{\OO}'_{\vec k}(X)$. Thus $T^1(X)$ is partitioned into disjoint subsets of generic vectors indexed by $\overline{\OO}'_{\vec k}(X)$, and summing their Liouville volumes yields the identity.

 \item In the four-holed sphere, the set of $\vec{1}$-maximal immersed pairs of pants of grading $\vec{1} = (1,1,1,1)$ coincides with the set of embedded pairs of pants, as explained in \cite{basmajian2025prime}. Thus, our identity reduces to the Luo–Tan identity in this case. However, on any other surface, the grading $\vec{1}$ index set of our identity and the Luo–Tan index set are generally quite different.

 \item When all grades are infinite, and the surface has at least one geodesic boundary component, then our identity recovers Bridgeman’s identity \cite{bridgeman2011orthospectra}.

 \item As the grade increases, $h$ decreases and the index set expands. 
 \item It may be possible to apply the method in \cite{doan2023measuring} to show monotonicity or convexity properties of the function $h$ in terms of the boundary length of an immersed pair of pants which could be exploited for applications. In the case of surfaces with only cusps, we obtain \cref{rem:Monotonicity} which gives us \cref{thm:Counting,Nondomination} as applications.

 \item It may be possible to relate our identity to the earlier identity in \cite{basmajian2025prime} via the boundary hitting function and its moments \cite{bridgeman2014moments}, but we haven't attempted this here. 
\end{enumerate}

\end{remark}
In order to describe concretely the function $h$ in \cref{thm:AbstractIdentity}, we recall the functions introduced by Bridgeman \cite{bridgeman2011orthospectra} (see also \cite{calegari2011bridgeman} for a different computation) and by Luo--Tan \cite{luo2014dilogarithm}. Consider the pre-immersed pair of pants $\widetilde{P}_{\vec k}(\mu)$ whose boundary components are $\widetilde{\alpha}^{k_{\alpha}}$, $\widetilde{\beta}^{k_{\beta}}$, and $\widetilde{\gamma}_{\vec{k}}(\mu)$ (see \cref{dfn:immersedPants}).

\begin{figure}[h!]
 \centering
 \begin{tikzpicture}[scale=2]
 \draw[thick] (-1,0) 
 arc[start angle=60, end angle=0, radius=1]; 

 \draw[thick] (0.5,{sqrt(3)/2}) 
 arc[start angle=60, end angle=0, radius=1]; 
 
 \draw[thick] (-0.5,-{sqrt(3)/2}) 
 arc[start angle=120, end angle=60, radius=1];

\draw[thick] (0.5,-{sqrt(3)/2}) arc[start angle=180, end angle=120, radius=1];
\node at (0.7,-0.5) {$\sigma_{\beta}$};

\draw[thick] (-1,0) arc[start angle=180, end angle=120, radius=1]; 
\node at (-0.7,-0.5) {$\sigma_{\alpha}$};

 \draw[thick] (1,0) 
 arc[start angle=240, end angle=180, radius=1]; 
 \draw[thick] (0.5,{sqrt(3)/2}) 
 arc[start angle=300, end angle=240, radius=1]; 

 \draw[thick] (0.5,-{sqrt(3)/2}) 
 arc[start angle=300, end angle=240, radius=1]; 
 \draw[dashed] (-0.5,{sqrt(3)/2}) 
 arc[start angle=0, end angle=-60, radius=1]; 
 \draw[thick]
 (-0.02,-0.73) to[out=90,in=260] (0,0.73);
 \draw[dashed]
 (0.02,-1) to[out=90,in=-80] (0,0.73);
 \node at (-0.2,0) {$\sigma_{\gamma}$};

 \node at (0.9,0.8) {$\widetilde{\beta}^{k_{\beta}}$};
 \node at (-0.8,0.8) {$\widetilde{\alpha}^{k_{\alpha}}$};
 \node at (0, -1.2) {$\widetilde{\gamma}_{\vec{k}}(\mu)$};
 \node at (0,0.9) {$\widetilde{\mu}$};
 \end{tikzpicture}
 \caption{Pre-immersed pair of pants \(\widetilde{P}_{\vec{k}}(\mu)\) (when $\alpha$ and $\beta$ are geodesics of finite grades).}
 \label{fig:Pants}
\end{figure}
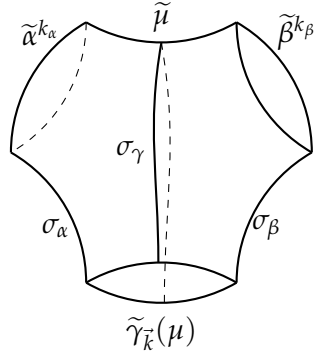

The Bridgeman function $Br$ associated to an orthogeodesic $\sigma$ is the Liouville measure of the set of unit tangent vectors generating geodesic arcs in the homotopy class of $\sigma$. On $\widetilde{P}_{\vec k}(\mu)$, we set 
\begin{equation}\label{equ:BrGamma}
 Br(\gamma_{\vec{k}}(\mu),\gamma_{\vec{k}}(\mu)) := 8\mathcal{L}\left(1/{\cosh^2\left(\frac{\ell\left(\sigma_{\gamma}\right)}{2}\right)}\right),
\end{equation} 
where $\sigma_{\gamma}$ is the simple orthogeodesic on $\widetilde{P}_{\vec k}(\mu)$ joining $\widetilde{\gamma}_{\vec{k}}(\mu)$ to itself. Similarly, \begin{equation}\label{equ:BrGamma2others}Br(\gamma_{\vec{k}}(\mu),\alpha^{k_{\alpha}}) := 8\mathcal{L}\left(1/{\cosh^2\left(\frac{\ell(\sigma_{\alpha})}{2}\right)}\right), \,\,\,\,\,\,\,\,\, Br(\gamma_{\vec{k}}(\mu),\beta^{k_{\beta}}) := 8\mathcal{L}\left(1/{\cosh^2\left(\frac{\ell(\sigma_{\beta})}{2}\right)}\right),
\end{equation}
where $\sigma_{\alpha}$ and $\sigma_{\beta}$ denote the two simple orthogeodesics (a.k.a. seams) on $\widetilde{P}_\mu$ connecting $\widetilde{\gamma}_{\vec{k}}(\mu)$ with $\widetilde{\alpha}^{k_{\alpha}}$ and $\widetilde{\beta}^{k_{\beta}}$ respectively (see \cref{fig:Pants}). Note that when $\widetilde{\alpha}^{k_{\alpha}}$ is a cone or a cusp, $Br(\gamma_{\vec{k}}(\mu),\alpha^{k_{\alpha}}) = 0$.

Following Luo–Tan \cite{luo2014dilogarithm}, when $\alpha^{k_{\alpha}}$ (respectively $\beta^{k_{\beta}}$) is a geodesic boundary component, the Lasso function $La$ associated to $\alpha^{k_{\alpha}}$ (respectively $\beta^{k_{\beta}}$) is defined to be the Liouville measure of the set of unit tangent vectors whose geodesic rays form lassos with footpoints on $\gamma_{\vec{k}}(\mu)$, which loop around $\alpha^{k_{\alpha}}$ (respectively $\beta^{k_{\beta}}$), and whose maximal geodesic extension is not homotopic (rel.\ boundary) to the simple orthogeodesic $\sigma_{\gamma}$. Explicitly, Luo–Tan \cite{luo2014dilogarithm} gives a closed-form expression for the Lasso function as follows:
\begin{equation}\label{equ:LassoAlpha}
La(\gamma_{\vec{k}}(\mu),\alpha^{k_{\alpha}})=2\Biggl(
 \mathcal{L}(y)
 -\mathcal{L}\Bigl(\frac{1 - x}{1 - x y}\Bigr)
 +\mathcal{L}\Bigl(\frac{1 - y}{1 - x y}\Bigr)
\Biggr),
\end{equation}
\begin{equation}\label{equ:LassoBeta}
La(\gamma_{\vec{k}}(\mu),\beta^{k_{\beta}})=2\Biggl(
 \mathcal{L}(z)
 -\mathcal{L}\Bigl(\frac{1 - x}{1 - x z}\Bigr)
 +\mathcal{L}\Bigl(\frac{1 - z}{1 - x z}\Bigr)
\Biggr),
\end{equation}
where $
x = e^{-\ell(\gamma_{\vec{k}}(\mu))},
y = \tanh^2\bigl(\tfrac{\ell(\sigma_{\alpha})}{2}\bigr),
$ and $z = \tanh^2\bigl(\tfrac{\ell(\sigma_{\beta})}{2}\bigr).$ When $\alpha$ is a cone point of angle $\theta(\alpha)$, then $\alpha^{k_{\alpha}}$ is a cone point of angle $\theta = k_{\alpha}\theta(\alpha)$. In \cref{sec:SmallConeAngle}, we establish the following (the case of $\beta^{k_{\beta}}$ is analogous).

\begin{proposition}\label{pro:LassoCuspCone}
The Lasso function $La(\gamma_{\vec{k}}(\mu),\alpha^{k_{\alpha}})$, in the case where $\alpha^{k_{\alpha}}$ is a cusp or a cone point of angle $\le \frac{\pi}{2}$, can be computed as follows:
\begin{enumerate}
 \item When $\alpha^{k_{\alpha}}$ is a cusp:
 \begin{equation}\label{eq:LassoCusp}
 La(\gamma_{\vec{k}}(\mu),\alpha^{k_{\alpha}})
 =4\,\mathcal{L}\left(\frac{1}{1+e^{\bar m_{\alpha}}}\right),
 \end{equation}
 where $\bar{m}_{\alpha}$ is the length of the truncated part (from $\gamma_{\vec{k}}(\mu)$ to the 1st natural collar of the cusp $\alpha^{k_{\alpha}}$) of the seam $\sigma_{\alpha}$.

 \item When $\alpha^{k_{\alpha}}$ is a cone point of angle $\theta=k_{\alpha}\theta(\alpha) \le \frac{\pi}{2}$:
 \begin{equation}\label{eq:LassoCone}
 \begin{aligned}
 La(\gamma_{\vec{k}}(\mu),\alpha^{k_{\alpha}})
 = 2\,\mathrm{Re}\left\{\mathcal{L}(z)-\mathcal{L}(1/z)\right\},
 \quad \text{where } 
 z : = \frac{\sinh m_{\alpha}+\tan(\theta/2)}{\sinh m_{\alpha} +i},
 \end{aligned}
 \end{equation}
 where $m_{\alpha}$ denotes the length of the seam $\sigma_{\alpha}$ from $\gamma_{\vec{k}}(\mu)$ to the cone point $\alpha^{k_{\alpha}}$.
\end{enumerate}
\end{proposition}

\begin{figure}[h!]
 \centering

 \begin{minipage}[t]{0.32\textwidth}
 \centering
 \begin{tikzpicture}[scale=0.35]
 \coordinate (C) at (-6,0); 
 \coordinate (T) at (5,1); 
 \coordinate (Tl) at (0,4); 
 \coordinate (B) at (5,-1); 
 \coordinate (G) at (0,-5); 
 \draw[thick]
 (C) to[out=0,in=0] (Tl);
 \draw[thick]
 (C) to[out=0,in=180] (Tl);
 \draw[thick]
 (Tl) to[out=5,in=180] (T);
 
 \draw[thick]
 (C) to[out=0,in=0] (G);
 \draw[thick]
 (C) to[out=0,in=-180] (G);
 \draw[thick]
 (G) to[out=-5,in=-180] (B);

 \draw[dashed, samples=120, smooth, domain=90:270, variable=\t]
 plot ({5 + 0.4*cos(\t)}, {0 + 1*sin(\t)});

 \draw[thick,samples=120, smooth, domain=-90:90, variable=\t]
 plot ({5 + 0.4*cos(\t)}, {0 + 1*sin(\t)});

 \node[left] at (0.5,3) {$\widetilde{\alpha}^{\infty}$};
 \node[above] at (6.5,-1) {$\widetilde{\beta}^{k_{\beta}}$};
 \node[above] at (-1,-2) {$\widetilde{\gamma}_{\vec{k}}(\mu)$};
 \node[above] at (2.4,0.8) {$\widetilde{\mu}$};
 \end{tikzpicture}
 \subcaption{Degeneration when $\alpha$ is a geodesic of grade $k_{\alpha}=\infty$ and $\beta$ is a geodesic of grade $k_{\beta}<\infty$.}
 \label{fig:GeodesicInfinite}
 \end{minipage}%
 \hfill
 \begin{minipage}[t]{0.32\textwidth}
 \centering
 \begin{tikzpicture}[scale=0.45]
 \coordinate (C) at (0,4); 
 \coordinate (T) at (5,1); 
 \coordinate (B) at (5,-1); 
 \coordinate (G) at (1,-3); 

 \draw[thick] (C) to[out=-90,in=180] (T);
 \draw[thick] (C) to[out=-90,in=180] (G);
 \draw[thick] (C) to[out=-90,in=0] (G);
 \draw[thick] (G) to[out=-5,in=180] (B);

 \draw[dashed, samples=120, smooth, domain=90:270, variable=\t]
 plot ({5 + 0.4*cos(\t)}, {0 + 1*sin(\t)});

 \draw[thick,samples=120, smooth, domain=-90:90, variable=\t]
 plot ({5 + 0.4*cos(\t)}, {0 + 1*sin(\t)});

 \node[left] at (C) {$\widetilde{\alpha}^{\infty}$};
 \node[right] at (5.6,0) {$\widetilde{\beta}^{k_{\beta}}$};
 \node[above] at (1,-3) {$\widetilde{\gamma}_{\vec{k}}(\mu)$};
 \node[above] at (2.4,1.3) {$\widetilde{\mu}$};
 \end{tikzpicture}
 \subcaption{Degeneration when $\alpha$ is a cusp of grade $k_{\alpha}=\infty$ and $\beta$ is a geodesic of grade $k_{\beta}<\infty$.}
 \label{fig:CuspInfinite}
 \end{minipage}%
 \hfill
 \begin{minipage}[t]{0.32\textwidth}
 \centering
 \begin{tikzpicture}[scale=0.35]
 \coordinate (C) at (-6,0); 
 \coordinate (T) at (5,4);
 \coordinate (Tu) at (1.5,3.5);
 \coordinate (Tl) at (-2,4);
 \coordinate (B) at (9,0);
 \coordinate (G) at (1.5,-1);
 \coordinate (Gu) at (1.5,1);

 \node[left] at (-1.3,3) {$\widetilde{\alpha}^{\infty}$};
 \node[right] at (4.6,3) {$\widetilde{\beta}^{\infty}$};
 \node[above] at (G) {$\widetilde{\gamma}_{\vec{k}}(\mu)$};
 \node[above] at (1.5,3.5) {$\widetilde{\mu}$};
 
 \draw[thick]
 (C) to[out=0,in=0] (Tl);
 \draw[thick]
 (Tl) to[out=0,in=180] (Tu);
 \draw[thick]
 (Tu) to[out=0,in=180] (T);
 \draw[thick]
 (T) to[out=0,in=180] (B);
 \draw[thick]
 (C) to[out=0,in=180] (Tl);

 \draw[dashed]
 (C) to[out=0,in=180] (Gu);
 \draw[dashed]
 (Gu) to[out=0,in=180] (B);
 \draw[thick]
 (C) to[out=0,in=-180] (G);
 \draw[thick]
 (G) to[out=0,in=-180] (B);
 \draw[thick]
 (T) to[out=180,in=-180] (B);
 \end{tikzpicture}
 \subcaption{Degeneration when $\alpha$ and $\beta$ are geodesics with $k_{\alpha}=k_{\beta}=\infty$.}
 \label{fig:GeodesicBothInfinite}
 \end{minipage}

 \caption{Degenerations of the pre-immersed pair of pants $\widetilde{P}_{\vec{k}}(\mu)$ in three limiting cases.}
 \label{fig:ThreeDegenerations}
\end{figure}
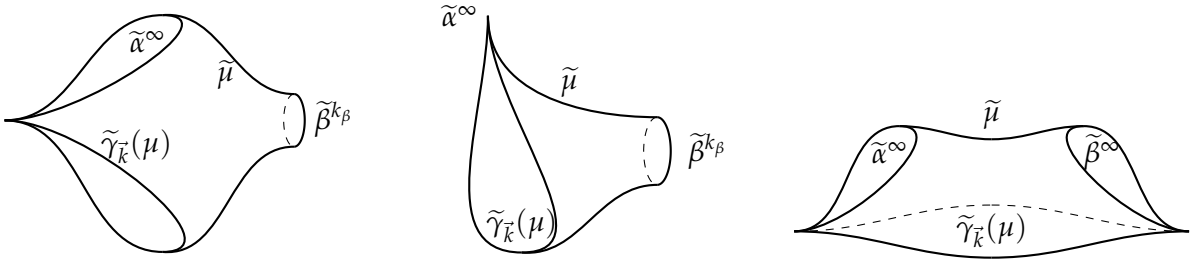

\begin{proposition}\label{prop:functionf}
Let $(X,\vec k)$ be an admissible pair such that, for every $i$ whose boundary element is a cone point, its angle satisfies $\theta_i \le \frac{\pi}{2k_i}$. Let $h$ be the function denoted in \cref{thm:AbstractIdentity}. Let $\mu$ be a $\vec{k}$-prime orthogeodesic connecting two boundary elements, say $\alpha$ and $\beta$, with grades $k_{\alpha}$ and $k_{\beta}$. If $k_{\alpha}$ and $k_{\beta}$ are finite (see \cref{fig:Pants}), then 
 $$h(P_{\vec{k}}(\mu)) = 4\pi^2-Br(\gamma_{\vec{k}}(\mu),\alpha^{k_{\alpha}})-Br(\gamma_{\vec{k}}(\mu),\beta^{k_{\beta}})-Br(\gamma_{\vec{k}}(\mu),\gamma_{\vec{k}}(\mu))$$
 $$-\, 4La(\gamma_{\vec{k}}(\mu),\alpha^{k_{\alpha}})-4La(\gamma_{\vec{k}}(\mu),\beta^{k_{\beta}}),$$ where all the terms on the RHS are given in \cref{equ:BrGamma,equ:BrGamma2others,equ:LassoAlpha,equ:LassoBeta,eq:LassoCusp,eq:LassoCone}. When at least one of the grades, say \(k_{\alpha}\), is infinite, the immersed pair of pants degenerates, and we obtain: 
\begin{enumerate}
 \item If $\alpha$ is a geodesic with grade $k_{\alpha} = \infty$ (see \cref{fig:GeodesicInfinite}), then 
 $$h(P_{\vec{k}}(\mu)) = 4\pi^2-Br(\gamma_{\vec{k}}(\mu),\beta^{k_{\beta}})-Br(\gamma_{\vec{k}}(\mu),\gamma_{\vec{k}}(\mu))-4La(\gamma_{\vec{k}}(\mu),\beta^{k_{\beta}}).$$

 \item If $\alpha$ is a cusp, and $k_{\alpha} = \infty$ (see \cref{fig:CuspInfinite}), then 
 $$h(P_{\vec{k}}(\mu)) = 0.$$

 \item If $k_{\alpha} = k_{\beta} = \infty$ (see \cref{fig:GeodesicBothInfinite}), then $h(P_{\vec{k}}(\mu))$ is the Bridgeman function, in particular 
 $$h(P_{\vec{k}}(\mu)) = Br(\alpha^{k_{\alpha}},\beta^{k_{\beta}}) = 8\mathcal{L}\left(1/{\cosh^2\left(\frac{\ell(\mu)}{2}\right)}\right).$$
\end{enumerate}
\end{proposition}

When the surface has only cusps, we have the following corollary. \begin{corollary}\label{cor:SimpleIdentity}
Let $(X,\vec k)$ be an admissible pair, where $X$ is a hyperbolic surface of genus $g$ with $n>0$ cusps. Then 
\begin{align*}
\sum_{\mu\in \overline{\OO}'_{\vec{k}}(X)} 
\Bigl(
\frac{\pi^2}{4}
-\mathcal L\Bigl(\frac{1}{e^{\bar m_{\alpha}}}\Bigr)
-\mathcal L\Bigl(\frac{1}{1+e^{\bar m_{\alpha}}}\Bigr)
\Bigr)
&=\frac{\pi^2}{4}(2g+n-2),
\end{align*}
where $\bar m_{\alpha}$ is as in \cref{eq:LassoCusp}.
\end{corollary}

By studying the monotonicity of the terms in these identities (see \cref{rem:Monotonicity}), we obtain the following two corollaries. They are inspired by \cite{doan2023measuring}, where a topological counting upper bound for embedded pairs of pants and non-domination of marked simple length spectra are derived from the Luo--Tan and Hu--Tan identities \cite{luo2014dilogarithm,hu2014new}. 

\begin{corollary}\label{thm:Counting}
Let $(X,\vec k)$ be an admissible pair, where $X$ is a hyperbolic surface of genus $g$ with $n>0$ cusps.
Let $\bar\mu$ be the truncation of a $\vec k$-prime orthogeodesic $\mu$ in the $\vec k$-concave core. Then, for $L>0$,
$$
\#\{\mu\in \overline{\OO}'_{\vec{k}}(X):\ell(\bar\mu)\le L\}
\le
\pi^2(2g+n-2)\frac{e^L}{L+2}.
$$
\end{corollary}

\begin{corollary}
\label{Nondomination}
Let $(X,\vec k)$ and $(Y,\vec k)$ be admissible pairs, where $X$ and $Y$ are hyperbolic surfaces of genus $g$ with $n>0$ cusps. Let $\bar\mu$ be the truncation of a $\vec k$-prime orthogeodesic $\mu$ in the $\vec k$-concave core. 
\begin{enumerate}
    \item If $
\ell_X(\bar\mu)\ge \ell_Y(\bar\mu)$ for all $\mu\in \overline{\OO}'_{\vec k}$,
then $X$ and $Y$ define the same point on the Teichm\"uller space, hence are isometric.

\item Let $\{\ell_X^{(j)}\}_{j\ge1}$ and $\{\ell_Y^{(j)}\}_{j\ge1}$ be the nondecreasing enumerations (with multiplicity) of the unmarked truncated $\vec k$-prime orthogeodesic length spectra. If $
\ell_X^{(j)}\ge \ell_Y^{(j)}$ for all $j$,
then $\ell_X^{(j)}=\ell_Y^{(j)}$ for all $j$.

\end{enumerate}
\end{corollary}


\begin{remark}
\begin{enumerate}
    \item We note that the identity in \cite[Thm.~1.5]{basmajian2025prime} after converting to $\ell(\bar\mu)$ also yields a uniform counting upper bound, but it is weaker than \cref{thm:Counting} (missing the $1/(L+2)$ factor). 
    \item \cref{Nondomination} also follows from \cite[Thm.~1.5]{basmajian2025prime} after rewriting in terms of $\ell(\bar\mu)$.
We note that, for surfaces of the same topological type with \emph{geodesic boundary}, Basmajian's and Bridgeman's identities imply that the full ortho-length spectrum admits no strict domination. See also \cite{masai2023systoles,quellec2024orthospectrum} for related rigidity results for the ortho-length spectrum.
\end{enumerate}

\end{remark}

\subsection{Organization}
After some preliminaries on the $k$-th natural collar and the model surface in \cref{sec:setup}, we prove the main abstract identity in \cref{sec:abstract}. We give the proof of the monotonicity and the two corollaries in \cref{sec:monotonicity}. In \cref{sec:ComputingTheMeasure}, we compute the measure appearing in the abstract identity and prove the simple identity in \cref{cor:SimpleIdentity}. In \cref{sec:SmallConeAngle}, we explicitly compute the closed forms of the Lasso functions associated to cusps and cones. In the final section, \cref{sec:Nongeneric}, we show that the set of generic vectors on hyperbolic cone surfaces has full Liouville measure.

\textbf{Acknowledgment.} Ara Basmajian is supported by PSC-CUNY Grant 68525-00 56 and a grant from the Simons foundation (SFI-MPS-TSM-00013865, A.B.). 
Nhat Minh Doan is funded by the Singapore National Research Foundation (NRF) under grant E-146-00-0029-01.
Hugo Parlier is supported by ANR-SNF Grant number 200021E\_238147 (SUGAR). Ser Peow Tan is supported by the National University of Singapore academic research grant A-8000989-00-00.

\section{Preliminaries and setup}\label{sec:setup}

In this section, we introduce the notation and describe the index sets for our identities, following the framework and conventions of Basmajian-Parlier-Tan \cite{basmajian2025prime} and reusing some of their results without proof.

\subsection{Notation}

We collect the various notations we use and the first place they appear in Table \ref{Table:notation}.


\begin{longtable}{| l | c | c |}
\caption{Definitions and notation}\label{Table:notation}\\
\hline\hline 
{\bf Definition} & {\bf Location} & {\bf Notation} \\
\hline\hline
\endfirsthead

\hline\hline 
{\bf Definition} & {\bf Location} & {\bf Notation} \\
\hline\hline
\endhead

\hline
\endfoot

\hline
\endlastfoot

Grading & \cref{eq:Grading} & $\vec{k} = (k_1,\hdots,k_n)$\\
\hline
Admissible pair & \cref{dfn:AdmissiblePair} & $(X,\vec k)$\\
\hline

$\vec k$-concave core & \cref{eq:Concave} & $V_{\vec{k}}(X)$ \\
\hline

Oriented $\vec k$-orthogeodesics on $X$ & Before \cref{dfn:immersedPants} & $\OO_{\vec k}(X)$\\
\hline

Unoriented $\vec k$-orthogeodesics on $X$ & Before \cref{dfn:immersedPants} & $\overline{\OO}_{\vec k}(X)$\\
\hline

$\vec k$-immersed pairs of pants associated to $\mu$ & \cref{dfn:immersedPants}& $P_{\vec k}(\mu)$\\
\hline

The closed geodesic associated to $\mu$ & \cref{eq:Gamma}& $\gamma_{\vec k}(\mu)$\\
\hline

Oriented $\vec{k}$-prime orthogeodesics & \cref{dfn:kPrimeOrthogeodesic}
& $\OO'_{\vec{k}}(X)$\\
\hline

Unoriented $\vec{k}$-prime orthogeodesics & \cref{dfn:kPrimeOrthogeodesic}
& $\overline \OO'_{\vec{k}}(X)$\\
\hline

\parbox[t]{0.55\linewidth}{%
Simple orthogeodesic on $\widetilde{P}_{\vec{k}}(\mu)$\\
joining $\widetilde{\gamma}_{\vec{k}}(\mu)$ to itself}
& \cref{equ:BrGamma}
& $\sigma_{\gamma}$ \\
\hline

\parbox[t]{0.55\linewidth}{%
Lasso function associated to $\alpha^{k_{\alpha}}$}
& \cref{equ:LassoAlpha}
& $La(\gamma_{\vec{k}}(\mu),\alpha^{k_{\alpha}})$ \\
\hline

$k$-th natural collar of $\delta$ &
\cref{dfn:kNaturalCollar} & $\CC_k(\delta)$ \\
\hline

$\vec{k}$-model surface & \cref{dfn:modelSurface} & $M_{\vec{k}}$ \\
\hline

Punctured cover of the model surface & \cref{dfn:punctureOrbifoldCover} & $\mathbb H_0$ \\
\hline

\end{longtable}

\subsection{The $k$-th natural collar of a boundary element}
Let $\beta$ be a boundary element on the hyperbolic cone surface $X$.
\begin{definition}\label{dfn:kNaturalCollar}{\cite[Definition 3.3]{basmajian2025prime}} The $k$-th natural collar $C_k(\beta)$ of $\beta$ is defined as \begin{enumerate}
\item the open horoball of boundary length $\frac{2}{k}$ around $\beta$ if $\beta$ is a cusp, \item the set of points at distance strictly less than $\arcsinh\left(1/\sinh(k\ell(\beta)/2)\right)$ from $\beta$ if $\beta$ is a simple closed geodesic. 
\item the set of points at distance strictly less than $\arccosh(1/\sin(k\theta/2))$ if $\beta$ is a cone point of angle $\theta$, where $\theta \le \pi/k$. \end{enumerate}
\end{definition}
In the cone case, $C_k(\beta)$ exists only when $\theta\le \pi/k$, and $C_k(\beta)=\varnothing$ when $\theta=\pi/k$. We have the following property.

\begin{proposition}[{\cite[Prop.~3.4]{basmajian2025prime}}]\label{pro:enterCollar}
Any geodesic path on $X$ that contains a subloop freely homotopic to the $k$-th power of $\beta$ necessarily intersects $C_k(\beta)$.
\end{proposition}

\subsection{Model surface and primality criterion}
Let $(X,\vec k)$ be an admissible pair. Let $\Sigma=\Sigma_{g,n}$ be a topological surface homeomorphic to $X$. Let $\overline{X}$ denote the metric completion of $X$ obtained by adding back all removed cone points.

\begin{definition}[Model surface $M_{\vec k}$]\label{dfn:modelSurface}
Given a grading vector $\vec k=(k_1,\dots,k_n)$, a \emph{$\vec k$–model surface} $M_{\vec k}$ is a hyperbolic/Euclidean orbifold on $\Sigma$ in which the $i$th boundary element corresponds to a deleted orbifold cone point of angle $\pi/k_i$ for $i\in\{1,\dots,n\}$. 
Let $\overline{M}_{\vec k}$ denote the metric completion obtained by adjoining all deleted cone points of positive angle (cusps are not adjoined).
\end{definition}

\begin{remark}[Existence]
A model surface $M_{\vec k}$ exists whenever the orbifold Euler characteristic is nonpositive (see \cite{Troyanov1991}):
$$
\chi_{\mathrm{orb}}(\Sigma_{g,n};\vec k)
=2(2-2g-n)+\sum_{i=1}^n \frac{1}{k_i}\le 0.
$$
If $\chi_{\mathrm{orb}}<0$ there are finite-area hyperbolic orbifold structures; if $\chi_{\mathrm{orb}}=0$ there are Euclidean orbifold structures. 
In particular, $M_{\vec k}$ exists for every admissible pair $(X,\vec k)$.
\end{remark}
\begin{remark}\label{rem:Retracing}
 Since every cone point of $\overline M_{\vec k}$ has even order, any geodesic passing through a cone point must retrace its path. 
\end{remark}

\begin{definition}[Collapsing map]\label{dfn:collapsingMap}
A \emph{collapsing map} is a surjective map
$$
f\colon \overline X \to \overline{M}_{\vec k}
$$
such that $f|_{\mathring{V}_{\vec k}(X)}:\mathring{V}_{\vec k}(X)\xrightarrow{\cong}M_{\vec k}$ is a marking-preserving homeomorphism and each component of $\overline X\setminus \mathring{V}_{\vec k}(X)$ is sent to the corresponding cone point in $\partial M_{\vec k} = \overline{M}_{\vec k} \setminus M_{\vec k}$. Here, $\mathring{V}_{\vec k}(X)$ denotes the interior of the concave core. 
\end{definition}

\begin{definition}[Punctured cover $\mathbb{H}_0$]\label{dfn:punctureOrbifoldCover} Let \[ \pi:\ \widetilde{\overline{M}_{\vec k}}\longrightarrow \overline{M}_{\vec k} \] be the universal orbifold covering map (so $\widetilde{\overline{M}_{\vec k}}=\mathbb H^2$ if $\chi_{\mathrm{orb}}<0$ and $\widetilde{\overline{M}_{\vec k}}=\mathbb E^2$ if $\chi_{\mathrm{orb}}=0$). Define the \emph{punctured cover} by \[ \mathbb{H}_0 \;:=\; \widetilde{\overline{M}_{\vec k}}\setminus \pi^{-1}\!\bigl(\{\text{cone points of }\overline{M}_{\vec k}\}\bigr). \] Restriction of $\pi$ gives a covering map \[ \pi:\ \mathbb{H}_0 \longrightarrow M_{\vec k} \] which is a \(2k_i\)-fold cover over each sufficiently small punctured disk around the orbifold cone point corresponding to the \(i\)-th boundary element of \(\Sigma\). \end{definition}


\begin{remark}\label{rem:homotopy-classes}
Under the collapsing map $f\colon \overline X\to\overline M_{\vec k}$, any geodesic arc $\alpha$ with endpoints in $\overline X\setminus \mathring{V}_{\vec{k}}(X)$ maps to either a cone point or a path $f(\alpha)$ joining the corresponding cone points; its homotopy class in $M_{\vec k}$ (rel.\ boundary) is denoted $[f(\alpha)]$.
 However, not every such homotopy class admits a properly immersed orthogeodesic representative in $M_{\vec k}$ - that is, an immersed geodesic segment whose interior lies entirely in $M_{\vec k}$ and whose endpoints are cone points. When this occurs, we say that $\alpha$ can be realised as a properly immersed orthogeodesic in $M_{\vec k}$.
\end{remark}
In fact, {\cite[Prop.~4.6]{basmajian2025prime}} states that:

\begin{proposition}[Primality criterion]\label{prop:prime-realization}
Let $(X,\vec k)$ be an admissible pair. An orthogeodesic in $\mathcal{O}_{\vec{k}}(X)$ is $\vec{k}$–prime if and only if it can be realized as a properly immersed orthogeodesic in $M_{\vec{k}}$.
\end{proposition}

\section{Proof of the abstract identity}\label{sec:abstract}

As noted in the introduction, \cref{thm:AbstractIdentity} follows from two inputs: 
(i) \cref{prop:uniqueness}, and 
(ii) the fact that the non-generic set $B_1\cup B_2\cup B_3$ (see \cref{dfn:NongenericVector}) has Liouville measure zero (see \cref{pro:WeakErgodicity}). 
In particular, \cref{prop:uniqueness} states that 
\begin{proposition}[\cref{prop:uniqueness}]
Let $(X,\vec k)$ be an admissible pair. Let $v\in T^1(X)$ be generic, and let $\omega_{[v]}$ be the geodesic path associated to $v$ as in \cref{dfn:genericVector}. 
Then there exists a unique $\vec k$–immersed pair of pants $P_{\vec k}(\mu)$, where $\mu$ is a $\vec k$–prime orthogeodesic, such that $\omega_{[v]}$ is isometrically immersed in $P_{\vec k}(\mu)$.
\end{proposition}

\begin{proof} 
The proof essentially follows the same argument as that of \cref{prop:prime-realization} in \cite[Prop.~4.6]{basmajian2025prime}, with additional concrete details. Let $M_{\vec{k}}$ be a model surface of $X$, and let $f:\overline X \to \overline{M}_{\vec{k}}$ be the collapsing map defined in \cref{dfn:collapsingMap}. For a generic vector $v\in T^1(X)$, denote by $\omega_{[v]}$ the geodesic path associated to $v$, oriented from the collar around a boundary element, say $\alpha$, to the collar around a boundary element, say $\alpha'$, where $\alpha$ and $\alpha'$ can be the same. We write $f(\omega_{[v]})$ for the image of $\omega_{[v]}$ under $f$. By construction, $\omega_{[v]}$ joins the two collars $\alpha$ and $\alpha'$ and meets no other collar in its interior. Thus $f(\omega_{[v]})$ is an immersed oriented arc in $M_{\vec{k}}$ connecting two cone points. Denote by $\widetilde{f(\omega_{[v]})}$ a lift of $f(\omega_{[v]})$ to the orbifold cover $\mathbb{H}_0$, in which orbifold points appear as punctures. 

\textbf{Existence.} We define the \emph{tightened representative} of $\widetilde{f(\omega_{[v]})}$ by homotoping $\widetilde{f(\omega_{[v]})}$ to "minimal" length within its homotopy class. As noted in \cref{rem:homotopy-classes}, not every such homotopy class admits a properly immersed minimal representative in $\mathbb{H}_0$. To overcome this obstacle, following \cite{basmajian2025prime}, we construct the representative infinitesimally as follows.

The tightened representative of $\widetilde{f(\omega_{[v]})}$ in $\mathbb{H}_0$, denoted by $[\widetilde{f(\omega_{[v]})}]$, can be represented as a sequence $\gamma_1\delta_1\gamma_2\delta_2\ldots\delta_{n-1}\gamma_n$ of directed geodesics $\gamma_i \subset \mathbb{H}_0$, between punctures $A_{i-1}$ and $A_i$ , $1\le i\le n$ and infinitesimal circular arcs $\delta_j, 1\le j \le n-1$, of turning angle $\theta_j$ about the intermediate punctures $A_j$ connecting the endpoint of $\gamma_j$ to the initial point of $\gamma_{j+1}$. Note that for $1 \le j \le n-1$, $|\theta_j| \ge \pi$, otherwise one can homotopic $\widetilde{f(\omega_{[v]})}$ away from that puncture to get a shorter representative. With this notation, we claim that
\begin{claim}\label{Claim:straighten}
If $n\ge 2$ then $|\theta_j| = \pi$, for all $1 \le j \le n-1$.
\end{claim}
\begin{figure}[h!]
 \centering
\begin{tikzpicture}[scale=1, rotate=0]
 
 \node[below] at (-4,0) {$\gamma_1$};
 \node[below] at (-2,0) {$\gamma_2$};
 \node[below] at (4,0) {$\gamma_{n}$};
 \node[above] at (-3,0.08) {$\delta_1$};
 \node[above] at (-1,0.08) {$\delta_2$};
 \node[above] at (3,0.08) {$\delta_{n-1}$};
 \node[below] at (-5,0) {$A_0$};
 \node[below] at (-3,0) {$A_1$};
 \node[below] at (-1,0) {$A_2$};
 \node[below] at (3,0) {$A_{n-1}$};
 \node[below] at (5,0) {$A_{n}$};

 \node[below] at (-7,0) {$A_{-1}$};
 \node[below] at (-6,0) {$\gamma_0$};
 \node[below] at (7,0) {$A_{n+1}$};
 \node[below] at (6,0) {$\gamma_{n+1}$};
\draw[dashed] (-7.7,0) to[out=0,in=180] (-7.1,0);
\draw[dashed] (7.1,0) to[out=0,in=180] (7.7,0);
 \draw[thin, mid arrow] (-2.9,0) to[out=0,in=180] (-1.1,0);
 \draw[thin, mid arrow] (-5,0) to[out=0,in=180] (-3.1,0);
 \draw[thin, mid arrow] (-7,0) to[out=0,in=180] (-5,0);
 \draw[thin, mid arrow] (5,0) to[out=0,in=180] (7,0);
 \draw[dashed] (-0.9,0) to[out=0,in=180] (2.9,0);
 \draw[thin, mid arrow] (3.1,0) to[out=0,in=180] (5,0);

 \draw[thin, mid arrow] (-5,0) to[out=20,in=160] (5,0);
 
 \node[above] at (0,1) {$\widetilde{f(\omega_{[v]})}$};

 \foreach \x in {-3,-1,3} {
 \draw[mid arrow]
 (\x,0) ++(-3pt,0)
 arc [start angle=180, end angle=0, radius=3pt];
 }

 \foreach \x in {-7,-5,-3,-1,3,5,7} {
 \draw[fill=white,draw=black] (\x,0) circle (1pt);
 }

\end{tikzpicture}
 \caption{Illustration of the lift $\widetilde{f(\omega_{[v]})}$ of $f(\omega_{[v]})$ to $\mathbb{H}_0$ and its tightened representative $\gamma_1\delta_1\gamma_2\delta_2\ldots\delta_{n-1}\gamma_n$.}

 \label{fig:lyingOneSide}
\end{figure}
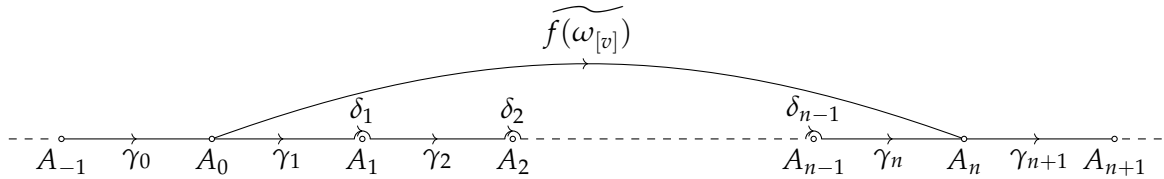
\begin{proof}[Proof of \cref{Claim:straighten}]
If $n\ge 2$ and $|\theta_j|>\pi$ for some $1 \le j \le n-1$, then, since $\mathbb{H}_0$ is a $2k_i$-fold cover of $M_{\vec{k}}$ over each small punctured disc around the $i$-th orbifold point, the projection of the tightened arc $[\widetilde{f(\omega_{[v]})}]$ to $M_{\vec{k}}$ is an arc, denoted by $[f(\omega_{[v]})]$, which contains a closed sub-loop that winds at least $k_j$ times around that orbifold point. Since $f(\omega_{[v]})$ is homotopic (rel.\ boundary) to its tightened representative $[f(\omega_{[v]})]$ in $M_{\vec{k}}$, the same sub-loop appears in $f(\omega_{[v]})$. Therefore, the preimage of $f(\omega_{[v]})$ under the marking-preserving homeomorphism $f|_{\mathring{V}_{\vec k}(X)}$ (defined in \cref{dfn:collapsingMap}) must also contain a closed sub-loop winding at least $k_j$ times around the corresponding boundary element of $X$. Here, the preimage under the map $f|_{\mathring{V}_{\vec k}(X)}$ of $f(\omega_{[v]})$ is the geodesic segment $\mathring{\omega}_{[v]}:=\omega_{[v]}\cap \mathring{V}_{\vec{k}}(X)$. By \cref{pro:enterCollar}, this forces $\mathring{\omega}_{[v]}$ to enter the $k_j$-th natural collar of that boundary element, contradicting the fact that $\mathring{\omega}_{[v]}\subset \mathring{V}_{\vec{k}}(X)$. Hence $|\theta_j|=\pi$ for every $1 \le j \le n-1$.
\end{proof}

Now we continue the proof of \cref{prop:uniqueness}. By \cref{prop:prime-realization}, one can find a $\vec{k}$-prime orthogeodesic $\mu$ in $X$ corresponding to the properly immersed oriented geodesic segment $\gamma_1$ connecting punctures $A_0$ and $A_1$ in $\mathbb{H}_0$. In particular, the arc $f(\mu)$ is homotopic to $\pi(\gamma_1)$ - the projection of $\gamma_1$ under the covering map $\pi:\mathbb{H}_0 \longrightarrow M_{\vec{k}}$, and $\mu$ runs from boundary element $\alpha$ to boundary element $\beta$, where $\alpha \subset f^{-1}(\pi(A_0))$ and $\beta \subset f^{-1}(\pi(A_1))$. By \cref{Claim:straighten}, all $\gamma_i$ lie on the same geodesic, say $g'\subset \overline{\mathbb{H}}_0=\mathbb H^2$ or $\mathbb E^2$ passing through $\{A_j\}_{j\in \mathbb{Z}}$. Thus by \cref{rem:Retracing}, all $\gamma_i$ project to $\pi(\gamma_1)$ in $M_{\vec{k}}$ with orientation alternating according to the parity of $i$. Hence, in $X$, the arc $\mathring{\omega}_{[v]}$ is homotopic (rel. boundary of $V_{\vec{k}}(X)$) to a (truncated) orthogeodesic $\eta$ of the form \begin{equation}\label{equ:formOfMu}
 \mu \beta^{\epsilon_1 k_{\beta}} \mu^{-1} \alpha^{\epsilon_2 k_{\alpha}}\mu\beta^{\epsilon_3 k_{\beta}} \ldots\mu^{s} ,
\end{equation}
in which $\mu$ appears $n$ times, $s:= (-1)^{n+1}$, and $\epsilon_j := \operatorname{sign}(\theta_j)$, for $1\le j \le n-1$. More precisely, the sign is chosen according to whether the orientation of the corresponding infinitesimal circular arc is counterclockwise (positive sign) or clockwise (negative sign). In fact, one can show that all~$\theta_j$ have the same sign (see \cref{fig:lyingOneSide} for the case $\theta_j=\pi$) by noting that if, for some~$j$, $\theta_j=\pi$ and either $\theta_{j-1}=-\pi$ or $\theta_{j+1}=-\pi$, then in~$X$ the arc $\mathring{\omega}_{[v]}$ would necessarily enter the collar. However, this stronger claim is not needed in our proof. 

Let $P_{\vec{k}}(\mu)$ be the immersed pair of pants associated to $\mu$. 
We note that, the natural collars arround $\alpha$ and $\beta$ lie inside $P_{\vec{k}}(\mu)$ as sets, and their preimages under the isometric immersion $\varphi:\widetilde{P}_{\vec{k}}(\mu)\to P_{\vec{k}}(\mu)$ are collars of $\widetilde{P}_{\vec{k}}(\mu)$, where $\widetilde{P}_{\vec{k}}(\mu)$ is the pre-immersed pair of pants of $P_{\vec{k}}(\mu)$ (see \cref{dfn:immersedPants}). Therefore, by the topological description of $\eta$ in \cref{equ:formOfMu}, every geodesic representative of its homotopy class (rel. boundary) - in particular $\mathring\omega_{[v]}$ and $\eta$ itself - is isometrically immersed in $P_{\vec k}(\mu)$. It follows that, $\omega_{[v]}$ is also isometrically immersed in $P_{\vec k}(\mu)$.

\textbf{Uniqueness.} 
Recall that the three curves 
$\alpha^{k_{\alpha}}, \beta^{k_{\beta}},$ and $\gamma_{\vec{k}}(\mu)$
form the boundary of the immersed pair of pants \(P_{\vec{k}}(\mu)\), where 
$$
\gamma_{\vec{k}}(\mu) \in [\alpha^{k_{\alpha}} * \mu * \beta^{k_{\beta}} * (\mu)^{-1}].
$$
On the punctured cover \(\mathbb{H}_0\), there are two distinct lifts of \(f(\gamma_{\vec{k}}(\mu))\), each tightening to a bi-infinite sequence
$$
\ldots \gamma_{-1}\,\delta_{-1}\,\gamma_0\,\delta_0\,\gamma_1\,\delta_1\ldots,
\quad\text{and}\quad
\ldots \gamma_{-1}\,\delta'_{-1}\,\gamma_0\,\delta'_0\,\gamma_1\,\delta'_1\ldots,
$$
where $
\delta_j$ and $\delta'_j$ are respectively infinitesimal circular arcs with turning angles $
 \theta_j = \pi$ and $\theta'_j = -\pi
$ for all $j\in\mathbb{Z}$.
These two lifts bound a bi-infinite strip, denoted by \(s_{\mu}\), which contains precisely the set of punctures \(\{A_j\}_{j\in\mathbb{Z}}\) on the geodesic \(g\subset \overline{\mathbb{H}}_0\). Moreover, since $\omega_{[v]}$ is isometrically immersed in $P_{\vec k}(\mu)$, the chosen lift \(\widetilde{f(\omega_{[v]})}\) of \(f(\omega_{[v]})\) (see \cref{fig:lyingOneSide}) lies entirely in the strip \(s_{\mu}\).

Now suppose \(\omega_{[v]}\) is also isometrically immersed in another \(\vec{k}\)-maximal pair of pants \(P_{\vec{k}}(\mu')\) with \(\mu'\neq \mu\in \overline{\mathcal{O}}'_{\vec{k}}(X)\). As before, choose a lift of \(f(\mu')\) in \(\mathbb{H}_0\) with one endpoint at \(A_0\) and the other at \(B\). This determines a bi-infinite strip \(s_{\mu'}\) bounded by two suitable lifts of $f(\gamma_{\vec{k}}(\mu'))$, which contains precisely the set of punctures on a geodesic \(g'\subset \overline{\mathbb{H}}_0\) passing through \(A_0\) and \(B\). Since \(\omega_{[v]}\) is isometrically immersed in \(P_{\vec{k}}(\mu')\), without loss of generality, we can assume that \(\widetilde{f(\omega_{[v]})}\subset s_{\mu'}\) as well. Then $s_{\mu}\cap s_{\mu'}$ contains the two distinct punctures $A_0$ and $A_n$, forcing $g=g'$ and hence $B=A_{1}$ or $B=A_{-1}$ due to \cref{prop:prime-realization}. Thus $f(\mu)$ is homotopic (rel. boundary) to $f(\mu')$ or $-f(\mu')$ in $M_{\vec{k}}$, and therefore $\mu$ and $\mu'$ are identical in $\overline{\mathcal{O}}'_{\vec{k}}(X)$, a contradiction.
\end{proof}

\section{Computing the measure}\label{sec:ComputingTheMeasure}
Let $\mu$ be a $\vec{k}$-prime orthogeodesic connecting two boundary elements, say $\alpha$ and $\beta$, with grades $k_{\alpha}$ and $k_{\beta}$, and let $c(\alpha)$ and $c(\beta)$ be their collars. Let $P_{\vec{k}}(\mu)$ be the immersed pair of pants associated to $\mu$ with three boundary components $\alpha^{k_{\alpha}},\beta^{k_{\beta}}$, and $\gamma_{\vec{k}}(\mu)$. In this section, we will prove \cref{prop:functionf} by computing the function $h(P_{\vec k}(\mu))$ defined in \cref{thm:AbstractIdentity}. 
 
We recall that the pre-immersed pair of pants $\widetilde{P}_{\vec{k}}(\mu)$ of $P_{\vec{k}}(\mu)$ is a hyperbolic pair of pants with boundary components \begin{equation}\label{eq:BoundaryToBoundary}
\widetilde{\alpha}^{k_{\alpha}}:=\varphi^{-1}(\alpha^{k_{\alpha}}), \widetilde{\beta}^{k_{\beta}}:=\varphi^{-1}(\beta^{k_{\beta}}),\widetilde{\gamma}_{\vec{k}}(\mu):=\varphi^{-1}(\gamma_{\vec{k}}(\mu))
\end{equation} (see \cref{fig:merged}) with first natural collars $c(\widetilde{\alpha}^{k_{\alpha}})$ and $c(\widetilde{\beta}^{k_{\beta}})$ of $\widetilde{\alpha}^{k_{\alpha}}$ and $\widetilde{\beta}^{k_{\beta}}$, respectively. The immersion $\varphi$ satisfies that \begin{equation}\label{eq:CollarToCollar}c(\widetilde{\alpha}^{k_{\alpha}})=\varphi^{-1}(c(\alpha)), \text{ and } c(\widetilde{\beta}^{k_{\beta}})=\varphi^{-1}(c(\beta)).
\end{equation}

We define $$V(\widetilde{P}_{\vec{k}}(\mu)):=\widetilde{P}_{\vec{k}}(\mu) \setminus (c(\widetilde{\alpha}^{k_{\alpha}}) \cup c(\widetilde{\beta}^{k_{\beta}}))$$ to be the concave core of the pre-immersed pair of pants $\widetilde{P}_{\vec{k}}(\mu)$.

Denote by $G(P_{\vec{k}}(\mu))$ the set of generic vectors $v \in T^1(X)$ with the associated geodesic arc $\omega_{[v]}$ isometrically immersed into $P_{\vec{k}}(\mu)$.

Similarly, we call a vector $u \in T^1(\widetilde{P}_{\vec{k}}(\mu))$ \emph{good} if its associated geodesic arc $\omega_{[u]}$ splits into three segments, the first and the last lie in the collar region $c(\widetilde{\alpha}^{k_{\alpha}}) \cup c(\widetilde{\beta}^{k_{\beta}})$, and the middle, denoted by $\mathring{\omega}_{[u]}$, lies in the concave core $V(\widetilde{P}_{\vec{k}}(\mu))$. Denote by $G(\widetilde{P}_{\vec{k}}(\mu))$ the set of all such good vectors $u \in T^1(\widetilde{P}_{\vec{k}}(\mu))$.

By the definition of $h$ in \cref{thm:AbstractIdentity}, we have that
$h(P_{\vec{k}}(\mu)) = \operatorname{Vol}(G(P_{\vec k}(\mu)))$. The following lemma lifts the computation of $\operatorname{Vol}(G(P_{\vec k}(\mu)))$ to $\operatorname{Vol}(G(\widetilde{P}_{\vec k}(\mu)))$ which is easier to deal with.
\begin{lemma}\label{lem:generic-lift} For every vector 
$v \in G\bigl(P_{\vec{k}}(\mu)\bigr)$, there exists a unique vector $u\in G(\widetilde{P}_{\vec{k}}(\mu))$ such that 
$v =(D\varphi)(u)$, where $D{\varphi}$ denotes the differential of $\varphi$.
\end{lemma}
\begin{proof}
Let $v\in G(P_{\vec{k}}(\mu))$ and suppose there are $u_1,u_2 \in G(\widetilde{P}_{\vec{k}}(\mu))$ with $v =(D\varphi)(u_1)=(D\varphi)(u_2)$. Then the associated geodesic arcs satisfy $\omega_{[v]}=\varphi(\omega_{[u_1]})=\varphi(\omega_{[u_2]})$. Thus, their subsegments lying in the concave core coincide:\begin{equation}\label{eq:UniquenessOfLift}
\mathring{\omega}_{[v]}=\varphi(\mathring{\omega}_{[u_1]})=\varphi(\mathring{\omega}_{[u_2]}). 
\end{equation}
Denote $H_{\alpha}:=\partial c(\widetilde{\alpha}^{k_{\alpha}})$ and $H_{\beta}:=\partial c(\widetilde{\beta}^{k_{\beta}})$ the two boundary components of the concave core $V(\widetilde{P}_{\vec k}(\mu))$ of the pre-immersed pair of pants. For each $i\in \{1,2\}$, without loss of generality, we may assume that $\mathring{\omega}_{[u_i]}$ is an oriented geodesic arc from $p_i \in H_{\alpha}$ to $q_i \in H_{\beta}$ making an angle $a_i$ with $H_{\alpha}$ and an angle $b_i$ with $H_{\beta}$. By \cref{eq:UniquenessOfLift}, \begin{equation}\label{eq:TwoPreThings}
 (\varphi(p_1),\varphi(q_1),a_1,b_1,\ell(\mathring{\omega}_{[u_1]}))=(\varphi(p_2),\varphi(q_2),a_2,b_2,\ell(\mathring{\omega}_{[u_2]})).
\end{equation}
Moreover, since $\varphi$ is $\pi_1$-injective and satisfies \cref{eq:BoundaryToBoundary,eq:CollarToCollar,eq:UniquenessOfLift}, the two geodesic arcs $\mathring{\omega}_{[u_1]}$ and $\mathring{\omega}_{[u_2]}$ must lie in the same homotopy class (rel. boundary of the concave core $V(\widetilde{P}_{\vec k}(\mu))$). \cref{equ:formOfMu} gives a topological description of their common image $\mathring{\omega}_{[v]}$ under $\varphi$. It follows that we can choose lifts $\widetilde H_{\alpha}$ of $H_{\alpha}$ and $\widetilde H_{\beta}$ of $H_{\beta}$ in the universal cover in $\mathbb{H}^2$ of $\widetilde{P}_{\vec k}(\mu)$ so that, for each $i \in \{1,2\}$, the corresponding lifted geodesic arc
$ \widetilde{\mathring\omega}_{[u_i]}$ of $ \mathring\omega_{[u_i]}$ runs from $\widetilde{p}_i \in \widetilde H_{\alpha}$ to $\widetilde{q}_i \in \widetilde H_{\beta}$. 
 
We now show that $\widetilde{p}_1=\widetilde{p}_2$ and $\widetilde{q}_1=\widetilde{q}_2$. There are two cases:
\begin{enumerate}
 \item If $ \widetilde{\mathring\omega}_{[u_1]} \cap \widetilde{\mathring\omega}_{[u_2]}\neq \varnothing$, then the two arcs together with $\widetilde H_{\alpha}$ and $\widetilde H_{\beta}$ bound two concave triangles in $\mathbb{H}^2$, each of which has total inner angles greater than $\pi$ due to \cref{eq:TwoPreThings}, a contradiction. 
 \item The second case when $ \widetilde{\mathring\omega}_{[u_1]} \cap \widetilde{\mathring\omega}_{[u_2]} = \varnothing$, then in $\mathbb{H}^2$, the four boundary curves and arcs bound a concave quadrilateral whose total interior angle equals to $2\pi$ due to \cref{eq:TwoPreThings}, again a contradiction. 
\end{enumerate}
Hence neither case is possible unless $\widetilde{p}_1=\widetilde{p}_2$ and $\widetilde{q}_1=\widetilde{q}_2$, as required. Consequently, $\omega_{[u_1]}$ and $\omega_{[u_2]}$ coincide. Both $u_1$ and $u_2$ are then tangent to the same geodesic arc.
Since \(D\varphi\) is injective on each tangent fiber, if \(u_1\neq u_2\) then they cannot lie over the same basepoint. If they lay over distinct basepoints, then the projected geodesic arc \(\omega_{[v]}\) would be tangent to itself at \(v\), a contradiction. Therefore $u_1 = u_2$.
\end{proof}
Before proving \cref{prop:functionf}, we need the following basic fact about hyperbolic geometry.
\begin{lemma}\label{lem:OutQ}
Let $Q\subset\mathbb H^2$ be the interior domain of an ideal quadrilateral with opposite sides $(a,c)$ and $(b,d)$, and let $\mu$ be the common perpendicular to $a$ and $c$. 
For $x\in\mu$, let $\sigma_1,\sigma_2$ be the geodesics through $x$ perpendicular to $b$ and $d$, and set $\phi:=\measuredangle(\sigma_1,\sigma_2)$.
If $\phi \le \pi/2$, then $x\notin Q$.
\end{lemma}

\begin{figure}[h]
\begin{center}
\begin{tikzpicture}[scale=2, line cap=round, line join=round]
 \pgfmathsetmacro{\rho}{0.6} 
 \pgfmathsetmacro{\x}{0.93} 
 \pgfmathsetmacro{\alpha}{\rho}
 \pgfmathsetmacro{\beta}{1/\rho}
 \pgfmathsetmacro{\h}{0.5*(\alpha+\beta)} 
 \pgfmathsetmacro{\r}{0.5*(\beta-\alpha)} 
 \pgfmathsetmacro{\y}{sqrt(1-\x*\x)} 

 \pgfmathsetmacro{\cOne}{(1+\alpha*\alpha)/(2*\x)}
 \pgfmathsetmacro{\cTwo}{(1+\beta*\beta)/(2*\x)}
 \pgfmathsetmacro{\Rone}{sqrt((\x-\cOne)*(\x-\cOne)+\y*\y)}
 \pgfmathsetmacro{\Rtwo}{sqrt((\x-\cTwo)*(\x-\cTwo)+\y*\y)}

 \pgfmathsetmacro{\Bx}{\alpha*\alpha/\cOne}
 \pgfmathsetmacro{\By}{\alpha*sqrt(1 - (\alpha*\alpha)/(\cOne*\cOne))}
 \pgfmathsetmacro{\Dx}{\beta*\beta/\cTwo}
 \pgfmathsetmacro{\Dy}{\beta*sqrt(1 - (\beta*\beta)/(\cTwo*\cTwo))}

 \pgfmathsetmacro{\phiOneStart}{atan2(\y, \x-\cOne)}
 \pgfmathsetmacro{\phiOneEnd}{atan2(\By, \Bx-\cOne)}
 \pgfmathsetmacro{\phiTwoStart}{atan2(\y, \x-\cTwo)}
 \pgfmathsetmacro{\phiTwoEnd}{atan2(\Dy, \Dx-\cTwo)}

 \draw[->] (-2.3,0) -- (2.3,0);
 \draw[thick] ({-\h+\r},0) arc (0:180:{\r});
 \draw[thick] ({\h+\r},0) arc (0:180:{\r});

 \draw[thick] ({\alpha},0) arc (0:180:{\alpha});
 \draw[thick] ({\beta},0) arc (0:180:{\beta});

 \draw[thick, dashed] (1,0) arc (0:180:1) node[pos=0.55, below left =1pt] {$\mu$};

 \draw[thick, dashed, teal] (0,\alpha) -- (0,\beta) node[pos=0.6, above left=2pt] {$\eta$};

 \draw[blue!70, thick] (\cOne,0) ++(\phiOneStart:\Rone) arc (\phiOneStart:\phiOneEnd:\Rone)
 node[pos=0.6, below= 0 pt] {$\sigma_1$};
 \draw[red!70, thick] (\cTwo,0) ++(\phiTwoStart:\Rtwo) arc (\phiTwoStart:\phiTwoEnd:\Rtwo)
 node[pos=0.55, right=2pt] {$\sigma_2$};

 \fill (0,1) circle (0.02) node[above right=1pt] {$y$};

 \fill (\x,\y) circle (0.02) node[ right=1pt] {$x$};

 \fill (\Bx,\By) circle (0.01);
 \fill (\Dx,\Dy) circle (0.01);

 \draw[black] (\Bx,\By) ++({atan2(\By,\Bx)}:0.06) -- ++({atan2(\By,\Bx)+90}:0.06);
 \draw[black] (\Bx,\By) ++({atan2(\By,\Bx)+90}:0.06) -- ++({atan2(\By,\Bx)+00}:0.06);

 \draw[black] (\Dx,\Dy) ++({atan2(\Dy,\Dx)+180}:0.06) -- ++({atan2(\Dy,\Dx)+90}:0.06);
 \draw[black] (\Dx,\Dy) ++({atan2(\Dy,\Dx)+90}:0.06) -- ++({atan2(\Dy,\Dx)+180}:0.06);

 \node at (-\h,0.4) {$a$};
 \node at ( 0.3+\h,0.3) {$c$};
 \node at (0,{\alpha-0.1}) {$b$};
 \node at (0,{\beta+0.1}) {$d$};

 \pgfmathsetmacro{\tOne}{atan2(\y,\x-\cOne)+90}
 \pgfmathsetmacro{\tTwo}{atan2(\y,\x-\cTwo)-90}
 \draw[thin] (\x,\y) ++(\tOne:0.08) arc (\tOne:\tTwo:0.08)
 node[midway, above=1pt] {$\phi$};
\end{tikzpicture}
\end{center}
 \caption{The ideal quadrilateral $Q$ and a point $x \notin Q$.}
 \label{fig:smallConeAngle}
\end{figure}
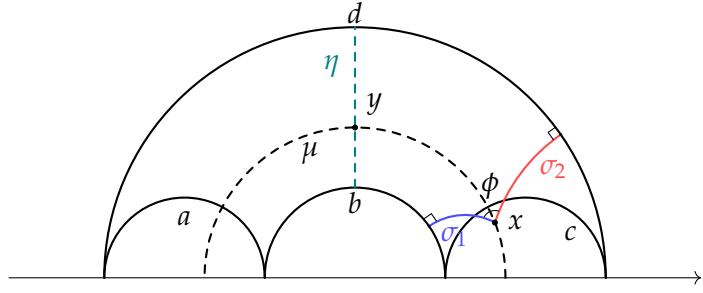

\begin{proof}
Normalize (see \cref{fig:smallConeAngle}) so that $\mu$ is the unit semicircle, $b,d$ are concentric semicircles, and $\eta$ is their common perpendicular, a diameter orthogonal to $\mu$. Let $y:=\mu\cap\eta$, $s:=d_{\mathbb H^2}(x,y)$, and $\delta:=d_{\mathbb H^2}(\mu,b)=d_{\mathbb H^2}(\mu,d)$. By symmetry, $\varphi:=\frac{\phi}{2}=\measuredangle(\mu,\sigma_1)$.

In the Lambert quadrilateral with sides $\mu,\eta,b,\sigma_1$,
\begin{equation}
\cot\varphi=\frac{\sinh s}{\coth\delta}
\qquad\text{\cite[Thm.~2.3.1(vi)]{buser2010geometry}}.
\end{equation}
Let $s':= d_{\mathbb H^2}(a,c)$. Symmetry gives $d_{\mathbb H^2}(y,a)=d_{\mathbb H^2}(y,c)=s'/2$, and in the Lambert quadrilateral with sides $b,\eta,\mu,c$, with one ideal vertex,
\begin{equation}
\sinh\left(\frac{s'}{2}\right)\sinh\delta=1 \qquad\text{\cite[Thm.~2.3.1(i)]{buser2010geometry}}.
\end{equation}

If $\phi \le \pi/2$ then $\varphi \le \pi/4$, so $\cot\varphi \ge 1$. From the second equation,
\[
\frac{\sinh(s'/2)}{\coth\delta}
=\frac{1}{\sinh\delta\,\coth\delta}
=\frac{1}{\cosh\delta}\le 1.
\]
Combining with the first gives $\sinh s \ge \sinh(s'/2)$, hence $s\ge s'/2$. Since $ \mu\cap Q=\{x\in\mu:\ s := d_{\mathbb H^2}(x,y)<s'/2\}$, we conclude $x\notin Q$.
\end{proof}

Now we prove the first part of \cref{prop:functionf}, namely the computation of the volume $h(P_{\vec{k}}(\mu))$ in the finite-grade case. The second part - concerning the degenerate cases where at least one grade is infinite - follows immediately (see \cref{fig:CuspInfinite,fig:GeodesicInfinite,fig:GeodesicBothInfinite}).

\begin{proposition} Let $(X,\vec k)$ be an admissible pair such that, for every $i$ whose boundary element is a cone point, its angle satisfies $\theta_i \le \frac{\pi}{2k_i}$. Let $\mu$ be a $\vec{k}$-prime orthogeodesic connecting two boundary elements, say $\alpha$ and $\beta$, with grades $k_{\alpha}$ and $k_{\beta}$. If $k_{\alpha}$ and $k_{\beta}$ are finite, then
 $$h(P_{\vec{k}}(\mu)) = \operatorname{Vol}(G(P_{\vec k}(\mu))) = 4\pi^2-Br(\gamma_{\vec{k}}(\mu),\alpha^{k_{\alpha}})-Br(\gamma_{\vec{k}}(\mu),\beta^{k_{\beta}})-Br(\gamma_{\vec{k}}(\mu),\gamma_{\vec{k}}(\mu))$$
\begin{equation}\label{eq:FormulaOfh}
 \hspace{1cm} -\,4La(\gamma_{\vec{k}}(\mu),\alpha^{k_{\alpha}})-4La(\gamma_{\vec{k}}(\mu),\beta^{k_{\beta}})
\end{equation} where all the terms on the RHS are given in \cref{equ:BrGamma,equ:BrGamma2others,equ:LassoAlpha,equ:LassoBeta}.
\end{proposition}
\begin{proof}
By \cref{lem:generic-lift}, the differential $D\varphi: T_{1}\bigl(\widetilde P_{\vec{k}}(\mu)\bigr) \longrightarrow T_{1}\bigl( P_{\vec{k}}(\mu)\bigr)$
restricts to a bijection from $G(\widetilde{P}_{\vec{k}}(\mu))$ to $G\bigl(P_{\vec{k}}(\mu)\bigr)$ and, moreover, preserves the Liouville volume on these measurable sets. Thus
$h(P_{\vec{k}}(\mu)) =\operatorname{Vol}(G(P_{\vec k}(\mu))) = \operatorname{Vol}(G(\widetilde{P}_{\vec k}(\mu))).$ Under \cref{dfn:AdmissiblePair}(1) of admissible pair $(X,\vec k)$, \cref{lem:OutQ} implies that the complement of $G(\widetilde{P}_{\vec k}(\mu))$ in $T^1(\widetilde{P}_{\vec k}(\mu))$ splits into 6 disjoint subsets $(A_i)_{i\in \{1,2,3,4,5,6\}}$ given by:
 \begin{enumerate}[label=($A_{\arabic*}$), ref=$A_{\arabic*}$]
 \item $v\in T_{1}(\widetilde P_{\vec k}(\mu))$ such that the geodesic arc 
 $\omega_{[v]}$ has one endpoint on $\widetilde\gamma_{\vec k}(\mu)$ 
 and the other on $\widetilde\alpha^{k_{\alpha}}$.
 
 \item $v\in T_{1}(\widetilde P_{\vec k}(\mu))$ such that the geodesic arc 
 $\omega_{[v]}$ has one endpoint on $\widetilde\gamma_{\vec k}(\mu)$ 
 and the other on $\widetilde\beta^{k_{\beta}}$.
 
 \item $v\in T_{1}(\widetilde P_{\vec k}(\mu))$ such that the maximal geodesic extension of 
 $\omega_{[v]}$ has both endpoints on $\widetilde\gamma_{\vec k}(\mu)$ 
 and is homotopic (rel.\ boundary) to the simple orthogeodesic 
 $\sigma_{\gamma}$ joining $\widetilde\gamma_{\vec k}(\mu)$ to itself 
 (see \cref{fig:Pants}).
 
 \item $v\in T_{1}(\widetilde P_{\vec k}(\mu))\setminus(A_{1}\cup A_{3})$ 
 such that $\omega_{[v]}$ has one endpoint on $\widetilde\gamma_{\vec k}(\mu)$ 
 and the other in the collar $c\bigl(\widetilde\alpha^{k_{\alpha}}\bigr)$.
 
 \item $v\in T_{1}(\widetilde P_{\vec k}(\mu))\setminus(A_{2}\cup A_{3})$ 
 such that $\omega_{[v]}$ has one endpoint on $\widetilde\gamma_{\vec k}(\mu)$ 
 and the other in the collar $c\bigl(\widetilde\beta^{k_{\beta}}\bigr)$.
 
 \item All remaining vectors in 
 $T_{1}\bigl(\widetilde P_{\vec k}(\mu)\bigr)\setminus\bigl(A_{1}\cup\cdots\cup A_{5} \cup G(\widetilde{P}_{\vec k}(\mu))\bigr)$.
\end{enumerate}
We observe that the subset \(A_{6}\) has Liouville measure zero, just as the non-generic set \(B_{1}\cup B_{2}\cup B_{3}\) introduced in \cref{dfn:NongenericVector}. 
We also observe that $v\in A_1$ if and only if the associated geodesic arc $\omega_{[v]}$ is homotopic to the simple orthogeodesic $\sigma_{\alpha}$ joining $\widetilde\gamma_{\vec k}(\mu)$ and $\widetilde\alpha^{k_{\alpha}}$ (see \cref{fig:Pants}), thus $\operatorname{Vol}(A_1)=Br(\gamma_{\vec{k}}(\mu),\alpha^{k_{\alpha}})$. Similarly, $\operatorname{Vol}(A_2)=Br(\gamma_{\vec{k}}(\mu),\beta^{k_{\beta}})$, and $\operatorname{Vol}(A_3)=Br(\gamma_{\vec{k}}(\mu),\gamma_{\vec{k}}(\mu))$. We also have that $\operatorname{Vol}(A_4)=4La(\gamma_{\vec{k}}(\mu),\alpha^{k_{\alpha}})$, since, for each vector $v\in A_4$, exactly one of $v,-v,r(v),-r(v)$ determines a lasso with footpoint on $\gamma_{\vec{k}}(\mu)$ which loops around $\alpha^{k_{\alpha}}$, where $r(v)$ denotes the reflection of $v$ across the three seams of the pair of pants. Similarly, $\operatorname{Vol}(A_5)=4La(\gamma_{\vec{k}}(\mu),\beta^{k_{\beta}})$. Thus \cref{eq:FormulaOfh} follows.
\end{proof}

\begin{proof}[Proof of \cref{cor:SimpleIdentity}]
Assume that $X$ has only cusps, so that $(X,\vec k)$ is an admissible pair with $X$ of genus $g$ and $n>0$ cusps. In this case the identity reduces to
\begin{equation}\label{eq:SimpleIdentity0}
\sum_{\mu\in \overline{\OO}'_{\vec{k}}(X)} 
\Bigl(
4\pi^2
-8\mathcal L\Bigl(\frac{1}{\cosh^2(\ell(\sigma_{\gamma})/2)}\Bigr)
-32\mathcal L\Bigl(\frac{1}{1+e^{\bar m_{\alpha}}}\Bigr)
\Bigr)
=4\pi^2(2g+n-2),
\end{equation}
where $\ell(\sigma_{\gamma})$, and $\bar m_{\alpha}$ are as in \cref{equ:BrGamma,eq:LassoCusp}.

Applying Abel’s five-term identity for the Rogers dilogarithm, specialized to $y=x$,
$$
2\mathcal L(x)=\mathcal L(x^2)+2\mathcal L\left(\frac{x}{1+x}\right),
$$
together with the hyperbolic relations in the pre-immersed pair of pants $P_{\vec k}(\mu)$ with two cusps (see \cite[Thm.~2.3.1]{buser2010geometry} and \cite[Lem.~7.3]{doan2025ortho}),
\begin{equation} \label{eq:TriFor}
\sinh\left(\frac{\ell\left(\gamma_{\vec k}(\mu)\right)}{4}\right)\sinh\left(\frac{\ell\left(\sigma_{\gamma}\right)}{2}\right)=1,\quad
e^{\bar m_{\alpha}}=\cosh\left(\frac{\ell\left(\sigma_{\gamma}\right)}{2}\right),\quad
e^{-\ell(\bar \mu)} = 1- e^{-2\bar m_{\alpha}},
\end{equation}
 the identity in \cref{eq:SimpleIdentity0} admits the following equivalent forms: 
\begin{align}
\sum_{\mu\in \overline{\OO}'_{\vec{k}}(X)} 
\Bigl(
\frac{\pi^2}{4}
-\mathcal L\Bigl(\frac{1}{e^{\bar m_{\alpha}}}\Bigr)
-\mathcal L\Bigl(\frac{1}{1+e^{\bar m_{\alpha}}}\Bigr)
\Bigr)
&=\frac{\pi^2}{4}(2g+n-2),\label{eq:graded-1}\\[0.4em]
\sum_{\mu\in \overline{\OO}'_{\vec{k}}(X)}
\Bigl(
\frac{\pi^2}{4}
-\mathcal L\Bigl(\frac{e^{\ell(\gamma_{\vec k}(\mu))/2}-1}{e^{\ell(\gamma_{\vec k}(\mu))/2}+1}\Bigr)
-\mathcal L\Bigl(\frac{1}{2}-\frac{1}{2e^{\ell(\gamma_{\vec k}(\mu))/2}}\Bigr)
\Bigr)
&=\frac{\pi^2}{4}(2g+n-2),\label{eq:graded-2}\\[0.4em]
\sum_{\mu\in \overline{\OO}'_{\vec{k}}(X)}
\Bigl(
\frac{\pi^2}{4}
-\mathcal L\bigl(\sqrt{1-e^{-\ell(\bar \mu)}}\bigr)
-\mathcal L\Bigl(\frac{\sqrt{1-e^{-\ell(\bar \mu)}}}{1+\sqrt{1-e^{-\ell(\bar \mu)}}}\Bigr)
\Bigr)
&=\frac{\pi^2}{4}(2g+n-2), \label{eq:graded-3}
\end{align}
where $\bar\mu$ is the truncation of the $\vec k$-prime orthogeodesic $\mu$ in the $\vec k$-concave core of $X$.
\end{proof}
\section{Monotonicity and applications}\label{sec:monotonicity}

We now show monotonicity properties of the functions in the identities of \cref{cor:SimpleIdentity} and deduce the corollaries stated in the introduction.
\begin{proposition}\label{rem:Monotonicity}
Each term
$
\frac{\pi^2}{4}
-\mathcal L\left(e^{-\bar m_{\alpha}}\right)
-\mathcal L\left((1+e^{\bar m_{\alpha}})^{-1}\right)$ in \cref{eq:graded-1} is a monotonic function of $\bar m_{\alpha}$, and takes values in $\bigl[0,\frac{\pi^2}{4}\bigr)$. Equivalently, each term in \cref{eq:graded-2,eq:graded-3} is monotonic in  $\ell(\gamma_{\vec k}(\mu))$ and $\ell(\bar\mu)$ respectively, where $\bar\mu$ is the truncation of the $\vec k$-prime orthogeodesic $\mu$ in the $\vec k$-concave core.
\end{proposition}
\begin{proof} The first part follows from the facts that $0 \le \bar m_{\alpha}<\infty$ and $\mathcal L$ is monotone. The second part follows from the relations $e^{-\ell(\bar\mu)} = 1 - e^{-2\bar m_{\alpha}}$ and $\ell(\gamma_{\vec k}(\mu))=4\,\operatorname{arcosh}(e^{\ell(\bar\mu)/2})$ (see \cref{eq:TriFor}).
\end{proof}
For $L>0$, define
\begin{equation}\label{eq:PhiF}
\Phi(L):=\frac{\pi^2}{4}
-\mathcal L\bigl(\sqrt{1-e^{-L}}\bigr)
-\mathcal L\Bigl(\frac{\sqrt{1-e^{-L}}}{1+\sqrt{1-e^{-L}}}\Bigr)    
\end{equation}

\begin{corollary}\label{thm:Counting1}
Let $(X,\vec k)$ be an admissible pair, where $X$ is a hyperbolic surface of genus $g$ with $n>0$ cusps. 
Then, for $L>0$,
$$
\#\{\mu\in \overline{\OO}'_{\vec{k}}(X):\ell(\bar\mu)\le L\}
\le
\frac{\frac{\pi^2}{4}(2g+n-2)}{\Phi(L)}
\le
\pi^2(2g+n-2)\frac{e^L}{L+2}.
$$
\end{corollary}

\begin{proof}
The first inequality follows from \cref{eq:graded-3} and the fact that $\Phi$ is strictly decreasing on $(0,\infty)$ (see \cref{rem:Monotonicity}). For the second, put $a:=\sqrt{1-e^{-L}}$. A standard Rogers--dilogarithm computation gives
\begin{equation}\label{eq:Phi-decomp}
\Phi(L)=\frac14\log\Bigl(\frac{1}{1-e^{-L}}\Bigr)\log\Bigl(\frac{1+a}{1-a}\Bigr)
+2\int_a^1\frac{-\log t}{1-t^2}\,dt.
\end{equation}
Using $-\log(1-x)\ge x$ and $-\log t\ge 1-t$,
\[
\frac14\log\Bigl(\frac{1}{1-e^{-L}}\Bigr)\ge \frac14 e^{-L},
\qquad
2\int_a^1\frac{-\log t}{1-t^2}\,dt
\ge 2\int_a^1\frac{1}{1+t}\,dt
=2\log\Bigl(\frac{2}{1+a}\Bigr)\ge 1-a.
\]
Since $1-a^2=e^{-L}$, we have $\log\bigl(\frac{1+a}{1-a}\bigr)=\log\bigl(\frac{(1+a)^2}{e^{-L}}\bigr)\ge L$, and
$\sqrt{1-x}\le 1-\tfrac{x}{2}$ with $x=e^{-L}$ gives $1-a\ge \tfrac12 e^{-L}$.
Substituting into \cref{eq:Phi-decomp} yields
\[
\Phi(L)\ge \frac{L}{4}e^{-L}+\frac12 e^{-L}=\frac{L+2}{4}e^{-L},
\]
which implies the second inequality.
\end{proof}

\begin{corollary}
\label{Nondomination1}
Let $(X,\vec k)$ and $(Y,\vec k)$ be admissible pairs, where $X$ and $Y$ are hyperbolic surfaces of genus $g$ with $n>0$ cusps.

\begin{enumerate}
    \item If $
\ell_X(\bar\mu)\ge \ell_Y(\bar\mu)$ for all $\mu\in \overline{\OO}'_{\vec k}$,
then $X$ and $Y$ define the same point on the Teichm\"uller space, hence are isometric.

\item Let $\{\ell_X^{(j)}\}_{j\ge1}$ and $\{\ell_Y^{(j)}\}_{j\ge1}$ be the nondecreasing enumerations (with multiplicity) of the unmarked truncated $\vec k$-prime orthogeodesic length spectra. If $
\ell_X^{(j)}\ge \ell_Y^{(j)}$ for all $j$,
then $\ell_X^{(j)}=\ell_Y^{(j)}$ for all $j$.

\end{enumerate}
\end{corollary}
\begin{proof}
For any admissible $(Z,\vec k)$, \cref{eq:graded-3} can be rewritten as \begin{equation}\label{eq:PhiSum}
\sum_{\mu\in \overline{\OO}'_{\vec k}} \Phi\bigl(\ell_Z(\bar\mu)\bigr)=\frac{\pi^2}{4}(2g+n-2).
\end{equation}

\begin{enumerate}
\item Assume $\ell_X(\bar\mu)\ge \ell_Y(\bar\mu)$ for all $\mu\in\overline{\OO}'_{\vec k}$.
Since $\Phi$ is decreasing, for every $\mu$ we have
\[
\Phi\bigl(\ell_X(\bar\mu)\bigr)\le \Phi\bigl(\ell_Y(\bar\mu)\bigr).
\]
Summing over $\mu$ and using \cref{eq:PhiSum} for $Z=X$ and $Z=Y$ shows that the two series are equal. Hence every termwise difference is zero, and
\[
\Phi\bigl(\ell_X(\bar\mu)\bigr)=\Phi\bigl(\ell_Y(\bar\mu)\bigr)
\qquad\text{for all }\mu.
\]
Because $\Phi$ is strictly decreasing, it follows that $\ell_X(\bar\mu)=\ell_Y(\bar\mu)$ for all $\mu\in\overline{\OO}'_{\vec k}$.
By \cref{prop:prime-realization}, there is a finite subset of $\overline{\OO}'_{\vec k}$ arising from a triangulation of the model surface; the associated truncated arc-length coordinates determine the hyperbolic structure. Therefore $X$ and $Y$ define the same point of Teichm\"uller space, and are isometric.

\item Let $\{\ell_X^{(j)}\}_{j\ge 1}$ and $\{\ell_Y^{(j)}\}_{j\ge 1}$ be the nondecreasing enumerations (with multiplicity) of the truncated unmarked $\vec k$-prime orthogeodesic lengths.
Assume $\ell_X^{(j)}\ge \ell_Y^{(j)}$ for all $j$.
Then $\Phi(\ell_X^{(j)})\le \Phi(\ell_Y^{(j)})$ for all $j$, and summing over $j$ gives
\[
\sum_{j\ge 1}\Phi(\ell_X^{(j)})\le \sum_{j\ge 1}\Phi(\ell_Y^{(j)}).
\]
Both sums equal the constant in \cref{eq:PhiSum} (they are the same series, just reordered), hence equality holds.
As above, this forces $\Phi(\ell_X^{(j)})=\Phi(\ell_Y^{(j)})$ for every $j$, and strict monotonicity of $\Phi$ implies
$\ell_X^{(j)}=\ell_Y^{(j)}$ for all $j$.
\end{enumerate}
\end{proof}
\section{Lasso functions on surfaces with small cone angles}\label{sec:SmallConeAngle}

\begin{figure}[h]
 \centering
\begin{tikzpicture}[scale=0.9]
 \pgfmathsetmacro{\x}{-3.0}
 \pgfmathsetmacro{\y}{-0.6}
 \pgfmathsetmacro{\c}{-1.2}
 \pgfmathsetmacro{\d}{ 2.8}

 \pgfmathsetmacro{\zOne}{0.5*(\x+\y)}
 \pgfmathsetmacro{\rOne}{0.5*(\y-\x)}
 \pgfmathsetmacro{\zTwo}{0.5*(\c+\d)}
 \pgfmathsetmacro{\rTwo}{0.5*(\d-\c)}

 \draw[->] (-4.5,0) -- (4.6,0);

 \draw[thick, name path=xy] (\y,0) arc[start angle=0, end angle=180, radius=\rOne] node[pos=.6, above] {$[x,y]$};
 \draw[thick, name path=cd] (\d,0) arc[start angle=0, end angle=180, radius=\rTwo] node[pos=.4, above] {$[c,d]$};

 \path[name intersections={of=xy and cd, by=q}];

 \coordinate (Z1) at (\zOne,0);
 \coordinate (Z2) at (\zTwo,0);

 \draw[dashed] (Z1) -- (q);
 \draw[dashed] (Z1) -- (Z2);
 \pic[draw, "$\phi$", angle radius=5pt, angle eccentricity=2] {angle = Z2--Z1--q};

 \fill (\x,0) circle (1pt) node[below] {$x$};
 \fill (\y,0) circle (1pt) node[below] {$y$};
 \fill (\c,0) circle (1pt) node[below] {$c$};
 \fill (\d,0) circle (1pt) node[below] {$d$};
 \fill (Z1) circle (1pt) node[below] {$z_1$};
 \fill (Z2) circle (1pt) node[below] {$z_2$};
 \fill (q) circle (1.2pt) node[above] {$q$};
\end{tikzpicture}

 \caption{The angle $\phi$.}
 \label{fig:AngleTheta}
\end{figure}
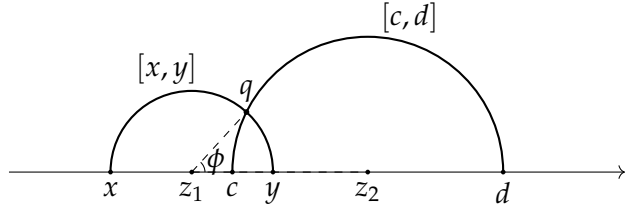

We recall a basic computation in the hyperbolic upper half-plane model $\mathbb{H}^2$.
\begin{lemma}[\cite{luo2014dilogarithm}, Lemma~4.2, Equation (12)]\label{lem:cotPhi}
Let $[x,y]$ and $[c,d]$ be complete geodesics in $\mathbb{H}^2$ with endpoints $x\neq y$, $c\neq d$ in $\partial \mathbb{H}^2 =\mathbb{R}\cup\{\infty\}$, intersecting at $q\in\mathbb{H}^2$. Let $z_1$ and $z_2$ be the Euclidean centers of the semicircles (or vertical lines, by limit) representing $[x,y]$ and $[c,d]$, respectively; when all endpoints are finite, $z_1=\tfrac{x+y}{2}$ and $z_2=\tfrac{c+d}{2}$. Set $\phi:=\measuredangle{(qz_1z_2)}$ (see \cref{fig:AngleTheta}). Then
\[
\cot^2\left(\frac{\phi}{2}\right)
=
\left|\frac{(x-c)(x-d)}{(y-c)(y-d)}\right|.
\]
\end{lemma}
 
\subsection{Lasso function associated to a cusp}
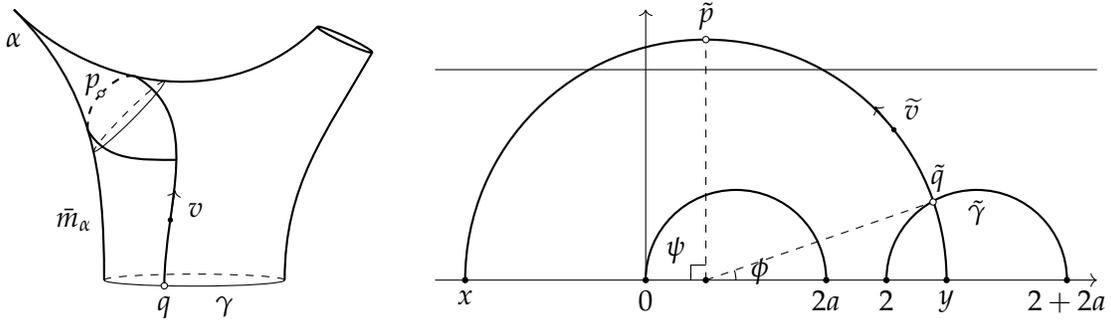
\begin{figure}[h!]
 \centering
 \begin{tikzpicture}[scale=0.4]
 

 
\coordinate (A) at (-6,6);
\coordinate (B) at (4.09,5.43);
\coordinate (C) at (-1,3.615);


 \draw[thin,dashed]
 (3,-3) 
 arc[x radius=3,y radius=0.2, start angle=0, end angle=180];
 \draw[thin]
 (3,-3) 
 arc[x radius=3,y radius=0.2, start angle=0, end angle=-180];

\draw[thick] (-3,-3) to[out=90,in=-45] (A);

\draw[thick] (A) to[out=-45,in=-135] (B);

\draw[thin,dashed, rotate around={45:(C)}]
 (C)
 arc[x radius=1.65, y radius=0.15, start angle=0, end angle=180];

\draw[thin, rotate around={45:(C)}]
 (C)
 arc[x radius=1.65, y radius=0.15, start angle=0, end angle=-180];

\draw[thick, rotate around={-25:(5,5)}] (5,5) ellipse (1 and 0.1);

\draw[thick] (3,-3) to[out=90,in=-120] (5.916,4.58);

\coordinate (D) at (-2,3.8);
\coordinate (E) at (-3.1,3.2);
\coordinate (F) at (-3.55,2);
\coordinate (G) at (-0.6,1);

\draw[thick] (-1,-3.19) to[out=90,in=-30] (D);

\draw[thick,dashed] (D) to[out=180,in=45] (E);

\draw[thick,dashed] (E) to[out=225,in=90] (F);

\draw[fill=white,draw=black] (E) circle (3pt);
\node at (-3.4,3.5) {$p$};

 \draw[fill=white,draw=black] (-1,-3.19) circle (3pt)
 node[below] {$q$};

\node at (-6,5) {$\alpha$};

\node at (-4,-1) {$\bar{m}_{\alpha}$};

 \draw[->,thin] (-0.8,-1) -- ++(81:1cm);
 \fill (-0.8,-1) circle (2.5pt);
\node[above right=4pt] at (-0.8,-1.5) {$v$};
 \draw[->,thin] (23.25,2) -- ++(129:1cm);
 \fill (23.25,2) circle (2.5pt);
 \node[above right=0pt] at (23.25,2) {$\widetilde {v}$};


\node[below] at (1,-3.19) {$\gamma$};

\draw[thick] (F) to[out=-60,in=-180] (G);


 \draw[->] (8,-3) -- (30,-3);

 \draw[->] (15,-3) -- (15,6);

 \draw[-] (8,4) -- (30,4);
 
 \draw[thick]
 (15,-3)
 arc[start angle=180, end angle=0, x radius=3, y radius=3];

 \draw[thick]
 (23,-3)
 arc[start angle=180, end angle=0, x radius=3, y radius=3];

 \draw[thick]
 (9,-3)
 arc[start angle=180, end angle=0, x radius=8, y radius=8];

 \fill (25,-3) circle (3pt) node[below] {$y$};

 \fill (9,-3) circle (3pt) node[below] {$x$};

 \draw[thin, dashed] (17,-3) -- (24.55,-0.4);

 \draw[thin] (18,-3) arc[start angle=0, end angle=19, radius=1cm];

 \node at (18.8,-2.7) {$\phi$};
 
 \draw[thin, dashed] (17,-3) -- (17,5);

 \draw[thin] (17,-2.5) -- (16.5,-2.5) -- (16.5,-3);
\node at (16,-2) {$\psi$};

 \draw[fill=white,draw=black] (17,5) circle (3pt) node[above] {$\tilde p$};
 
 \draw[fill=white,draw=black] (24.55,-0.4) circle (3pt); 
 
 \node[above] at (24.7,-0.4) {$\tilde q$};

 \node[below] at (26,0) {$\tilde \gamma$};


 \fill (15,-3) circle (3pt) node[below] {$0$};

 \fill (17,-3) circle (3pt);
 
 \fill (21,-3) circle (3pt) node[below] {$2a$};

 \fill (23,-3) circle (3pt) node[below] {$2$};

 \fill (29,-3) circle (3pt) node[below] {$2+2a$};

\end{tikzpicture}
 \caption{A lasso about a cusp and its lift to $\mathbb{H}^2$.}

 \label{fig:LassoCusp}
\end{figure}

We can now calculate the Lasso function associated to a cusp. Let $\alpha$ be the cusp and $\gamma$ be the other boundary of the pair of pants on which the foot of the lasso lies. Let $\bar{m}_{\alpha}$ be the length of the seam from $\gamma$ to the natural collar of the cusp $\alpha$, see \cref{fig:LassoCusp}. We can lift the picture to the universal cover, see \cref{fig:LassoCusp}, where the cusp is lifted to infinity and $[0,2a]$ and $[2,2+2a]$ are lifts of $\gamma$, and $a=e^{-\bar{m}_{\alpha}}$.

If $v$ is a unit tangent vector which exponentiates to a lasso we are interested in, let $[x,y]$, where $x<0$, $y \in (2,2+2a)$ be a geodesic in $\mathbb{H}^2$, which is the relevant lift of the extension of the lasso to a complete geodesic in $\mathbb{H}^2$, and $\tilde p, \tilde q \in [x,y]$, where $\tilde p$ is the point closest to the cusp at $\infty$, so the highest point on the semicircle $[x,y]$ and $\tilde q$ is the intersection of $[x,y]$ with $[2,2+2a]$. Then $v$ projects to a point between $\tilde p$ and $\tilde q$.

\begin{lemma}\label{lem:LassoCusp}
The measure of the set of unit tangent vectors giving rise to lassos with the foot on $\gamma$ and looping counterclockwise around the cusp $\alpha$ (see \cref{fig:LassoCusp}) is
\[
La(\gamma,\alpha)
= \int_{-\infty}^{0}\int_{2}^{2+2a}\frac{2\,d_{\mathbb H^2}(\tilde p,\tilde q)}{(y-x)^2}\,dy\,dx
= \int_{-\infty}^{0}\int_{2}^{2+2a}\frac{2\ln\cot(\phi/2)}{(y-x)^2}\,dy\,dx
\]
\[
= \int_{-\infty}^{0}\int_{2}^{2+2a}\frac{\ln\left|\frac{(2+2a-x)(2-x)}{(y-2)(2+2a-y)}\right|}{(y-x)^2}\,dy\,dx
= 4\,\mathcal L\left(\frac{a}{1+a}\right) = 4\,\mathcal L\left(\frac{1}{e^{\bar m_{\alpha}}+1}\right),
\]
where $a:=e^{-\bar m_{\alpha}}$ and $\mathcal L$ is the Rogers dilogarithm.
\end{lemma}

\begin{proof}
The hyperbolic distance between $e^{i\phi}$ and $-e^{-i\psi}$ in $\mathbb H^2$
(with $\phi,\psi\in(0,\pi/2]$) is
\[
d_{\mathbb H^2}\bigl(e^{i\phi},-e^{-i\psi}\bigr)
= \ln\cot(\phi/2)+\ln\cot(\psi/2).
\]
In our case $\psi=\pi/2$, so $\cot(\pi/4)=1$ and hence
\[
d_{\mathbb H^2}(\tilde p,\tilde q)
= \ln\cot(\phi/2)
= \frac{1}{2}\,\ln\left(\frac{(2+2a-x)(2-x)}{(y-2)(2+2a-y)}\right), \text{ by \cref{lem:cotPhi}.}
\]

For any $c\in\mathbb R$,
\[
\int \frac{\ln|y-c|}{(y-x)^2}\,dy
= -\frac{\ln|y-c|}{\,y-x\,}+\frac{\ln|y-x|}{x-c}-\frac{\ln|y-c|}{x-c}+\mathrm{const},
\]
and
\[
\int \frac{1}{(y-x)^2}\,dy = -\frac{1}{\,y-x\,}+\mathrm{const}.
\]
Hence, for $c\neq d$,
\begin{align*}
\int \frac{\ln\left|\frac{(x-c)(x-d)}{(y-c)(y-d)}\right|}{(y-x)^2}\,dy
&= \Bigl(\frac{1}{y-x}+\frac{1}{x-c}\Bigr)\ln|y-c|-\left(\frac{1}{x-c}+\frac{1}{x-d}\right)\ln|y-x| \\
&\quad +\Bigl(\frac{1}{y-x}+\frac{1}{x-d}\Bigr)\ln|y-d|
 - \frac{\ln|(x-c)(x-d)|}{\,y-x\,}+\mathrm{const},
\end{align*}
and evaluating from $y=c$ to $y=d$ yields
\[
\int_{c}^{d} \frac{\ln\left|\frac{(x-c)(x-d)}{(y-c)(y-d)}\right|}{(y-x)^2}\,dy
= \frac{2}{x-d}\ln\left|\frac{x-c}{d-c}\right|
 - \frac{2}{x-c}\ln\left|\frac{x-d}{c-d}\right|
= 2\,\frac{d}{dx}\,J(x;d,c),
\]
where \(J(x;d,c):=2\,\mathcal L\bigl(\frac{x-c}{d-c}\bigr)\).
Taking \(c=2\), \(d=2+2a\) and integrating in $x$ from $-\infty$ to $0$ gives
\[
\int_{-\infty}^{0}\int_{2}^{2+2a}
\frac{\ln\left|\frac{(2+2a-x)(2-x)}{(y-2)(2+2a-y)}\right|}{(y-x)^2}\,dy\,dx
= 2\int_{-\infty}^{0}\frac{d}{dx}J(x;2+2a,2)\,dx
\]
\[
= 4\,\mathcal{L}\left(-\frac{1}{a}\right)
 - 4\lim_{x\to -\infty}\mathcal{L}\left(\frac{x-2}{2a}\right)
= -4\,\mathcal L\left(\frac{1}{1+a}\right) + 4\,\mathcal{L}(1) = 4\,\mathcal L\left(\frac{a}{1+a}\right),
\]
by Euler's reflection identity $\mathcal L(t)+\mathcal L(1-t)= \mathcal L(1) = \pi^2/6$ and Landau's identity \(\mathcal L(-1/t)=-\mathcal L\bigl(\frac{1}{1+t}\bigr)\) for \(t>0\).
\end{proof}
\begin{remark}
Another way to obtain the Lasso function associated to a cusp is as the limit of the Lasso function associated to a geodesic boundary component in \cref{equ:LassoAlpha}, as the length of the geodesic boundary tends to $0$.

\end{remark}

\subsection{Lasso function associated to a cone point}\label{subsec:SmallConeAngle}
\begin{figure}[h!]
 \centering
 \begin{tikzpicture}[scale=0.4]
 

 
\coordinate (A) at (-4.2,4.6);
\coordinate (B) at (4.09,5.43);
\coordinate (C) at (-1,3.615);


 \draw[thin,dashed]
 (3,-3) 
 arc[x radius=3,y radius=0.2, start angle=0, end angle=180];
 \draw[thin]
 (3,-3) 
 arc[x radius=3,y radius=0.2, start angle=0, end angle=-180];

\draw[thick] (-3,-3) to[out=90,in=-75] (A);

\draw[thick] (A) to[out=-25,in=-135] (B);

\draw[thin,dashed, rotate around={45:(C)}]
 (C)
 arc[x radius=1.7, y radius=0.15, start angle=0, end angle=180];

\draw[thin, rotate around={45:(C)}]
 (C)
 arc[x radius=1.7, y radius=0.15, start angle=0, end angle=-180];

\draw[thick, rotate around={-25:(5,5)}] (5,5) ellipse (1 and 0.1);

\draw[thick] (3,-3) to[out=90,in=-120] (5.916,4.58);

\coordinate (D) at (-2,3.8);
\coordinate (E) at (-3.1,3.2);
\coordinate (F) at (-3.55,2);
\coordinate (G) at (-0.6,1);

\draw[thick] (-1,-3.19) to[out=90,in=-30] (D);

\draw[thick,dashed] (D) to[out=180,in=45] (E);

\draw[thick,dashed] (E) to[out=225,in=90] (F);

\draw[fill=white,draw=black] (E) circle (3pt);
\node at (-3.4,3.5) {$p$};


 \draw[fill=white,draw=black] (-1,-3.19) circle (3pt)
 node[below] {$q$};

\node at (-4.8,5) {$\alpha$};
\node at (-4.5,1) {$m_{\alpha}$};


\node[below] at (1,-3.19) {$\gamma$};

\draw[thick] (F) to[out=-60,in=-180] (G);


 \draw[->] (8,-3) -- (30,-3);

 \draw[->] (18,-3) -- (18,6);

 
 \draw[thick]
 (13.5,-3)
 arc[start angle=180, end angle=0, x radius=1.75, y radius=1.75];

 \draw[thin, dashed]
 (16.6,-3)
 arc[start angle=180, end angle=80, x radius=11.15, y radius=11.15];
 
 \node at (27,7) {$[-1/k, k]$};

 \draw[thick]
 (19,-3)
 arc[start angle=180, end angle=0, x radius=1.75, y radius=1.75];

 \draw[thick]
 (11,-3)
 arc[start angle=180, end angle=0, x radius=5, y radius=5];
 
 \draw (18,3.3) circle (2.3); 
 \draw (18,2.8) circle (1.1); 


 \fill (18,2.4) circle (3pt) node[right] {$i$};

 \fill (21,-3) circle (3pt) node[below] {$y$};

 \fill (11,-3) circle (3pt) node[below] {$x$};

 \draw[fill=white,draw=black] (17.62,1.72) circle (3pt) node[below] {$\tilde p$};

 \draw[fill=white,draw=black] (20.68,-1.26) circle (3pt) node[below right] {$\tilde q$};

 \draw[->,thin] (-0.8,-1) -- ++(81:1cm);
 \fill (-0.8,-1) circle (2.5pt);
\node[above right=4pt] at (-0.8,-1.5) {$ v$};
 \draw[->,thin] (20.03,0) -- ++(129:1cm);
 \fill (20,0) circle (2.5pt);
 \node[above right=0pt] at (20,0) {$\widetilde { v}$};

 \node[below] at (22.2,-0.5) {$\tilde \gamma$};


 \fill (13.5,-3) circle (3pt) node[below] {$-b$};

 \fill (17,-3) circle (3pt) node[below] {$-a$};

 \fill (19,-3) circle (3pt) node[below] {$a$};

 \fill (22.5,-3) circle (3pt) node[below] {$b$};

 \fill (18,-3) circle (3pt) node[below] {$0$};

\end{tikzpicture}
 \caption{A lasso about a cone point and its lift to $\mathbb{H}^2$.}

 \label{fig:LassoCone}
\end{figure}
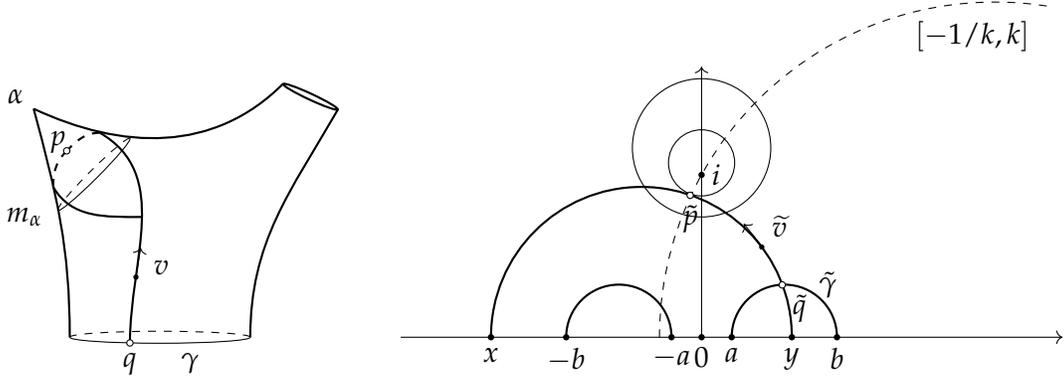

We next calculate the Lasso function associated to a cone point with cone angle $\theta$. In order for the formula for the various measures to hold over all elements of the index set, we will require $\theta \leq \pi/2$ (see \cref{lem:OutQ}). Again, we let $\alpha$ be the cone point, $\gamma$ be the other component of the boundary of the pair of pants on which the foot of the lasso lies, $m_{\alpha}$ be the distance from $\alpha$ to $\gamma$. Denote by $La(m_{\alpha}, \theta)$ the measure of the set of unit tangent vectors corresponding to the lasso. Take a lift to $\mathbb{H}^2$ so that $\alpha$ lifts to $i=\sqrt{-1}$, and two of the lifts for $\gamma$ are $[a,b]$ and $[-b,-a]$, where $0<a<b$. Here $a,b$ are functions of $\theta, m_{\alpha}$, and if $\theta\leq \pi/2$, then by \cref{lem:OutQ}, $b<1$. See \cref{fig:LassoCone}.

 The relevant lift of the extension of the lasso to $\mathbb{H}^2$ (via the developing map) is the geodesic $[x,y]$ with
$$a<y<b,\qquad -\frac{1}{y}<x<-b.$$
The bound $-\frac{1}{y}<x$ ensures that the geodesic $[x,y]$ lies below the geodesic $[-1/y,y]$ passing through $i=\sqrt{-1}$. The foot of $v$ lies in the part between the points $\tilde p$ and $\tilde q$ where $\tilde q=[x,y]\cap[a,b]$, and $\tilde p$ is the point on $[x,y]$ closest to $i$. In particular, $\tilde p$ is the point at the intersection of $[x,y]$ with the geodesic $[-1/k,k]$, where $k$ satisfies the harmonic cross ratio (see \cite[§7.24]{beardon2012geometry})$$[x,y,-1/k,k]=\frac{(x+1/k)(y-k)}{(x-k)(y+1/k)}=-1.$$ Therefore, $k$ can be expressed in terms of $x$ and $y$.
The relation is given by 
\begin{equation}\label{eq:quadratic in k}
 \frac{1}{k}-k=\frac{2-2xy}{x+y}.
\end{equation}

We have:
\begin{lemma}\label{lem:LassoCone1}
	The measure of the set of unit tangent vectors giving rise to lassos with the foot on $\gamma$ and looping counterclockwise around the cone point $\alpha$ (see \cref{fig:LassoCone}) is given by 
	\begin{align*}
		La(\gamma,\alpha)=&\int_{a}^{b}\int_{-1/y}^{-b}\frac{2d_{\mathbb H^2}(\tilde p,\tilde q)}{(y-x)^2} \,dxdy
	=\int_{a}^{b}\int_{-1/y}^{-b}
\frac{\ln\left|\frac{(x-a)(x-b)(1+y^2)}{(y-a)(y-b)(1+x^2)}\right|}{(y-x)^2}\,dx\,dy,
	\end{align*}
	where $a,b$ can be expressed in terms of the parameters $\theta$ and $m_{\alpha}$. In particular,
 $$a=\frac{\cosh m_{\alpha}\cos (\theta/2)-1}{\sinh m_{\alpha}+\cosh m_{\alpha}\sin (\theta/2)},\qquad
b=\frac{\cosh m_{\alpha}\cos (\theta/2)+1}{\sinh m_{\alpha}+\cosh m_{\alpha} \sin (\theta/2)}.$$
\end{lemma}
\begin{proof}
 Let $\phi$ and $\psi$ be the angles subtended by $\tilde q$ and $\tilde p$ from $(x+y)/2$ with the positive real line, then $$d_{\mathbb H^2}(\tilde p,\tilde q)=\ln \cot(\phi/2)-\ln \cot(\psi/2) = \frac{1}{2}\ln \left|\frac{(x-a)(x-b)(y+1/k)(y-k)}{(y-a)(y-b)(x+1/k)(x-k)}\right|,$$ by \cref{lem:cotPhi}, where $k$ is given in \cref{eq:quadratic in k}. Using 
\cref{eq:quadratic in k},
\begin{displaymath}
\frac{(y+1/k)(y-k)}{(x+1/k)(x-k)}=\frac{y^2-1+(1/k-k)y}{x^2-1+(1/k-k)x}=\frac{(x+y)(y^2-1)+(2-2xy)y}{(x+y)(x^2-1)+(2-2xy)x}=-\frac{y^2+1}{x^2+1}.
\end{displaymath}
Substituting into the previous display gives the stated integral formula. The formulas for $a$ and $b$ in terms of $m_{\alpha}$ and $\theta$ involve standard calculations in hyperbolic geometry and will be omitted.
\end{proof}
Thus, we are left to compute
\begin{equation}\label{ }
I(a,b):=\int_{a}^{b}\int_{-1/y}^{-b}\frac{\ln \left|\frac{(x-a)(x-b)(y^2+1)}{(y-a)(y-b)(x^2+1)}\right|}{(y-x)^2} \,dxdy,
\end{equation}
where $0<a<b<1$. Recall that the Rogers dilogarithm $\mathcal{L}(\,\cdot\,)$ is defined holomorphically on $\mathbb{C}\setminus\big((-\infty,0]\cup[1,\infty)\big)$ by
$$
\mathcal{L}(z) := \Li_2(z)+\frac12 \Log(z)\,\Log(1-z) = -\frac12\int_{0}^{z}\left(\frac{\Log(1-t)}{t}+\frac{\Log(t)}{1-t}\right)dt,
$$
where $\Li_2$ is the principal branch of the dilogarithm (cut along $[1,\infty)$) and $\Log$ is the principal logarithm (cut along $(-\infty,0]$). The integral is taken along any path in $\mathbb{C}\setminus\big((-\infty,0]\cup[1,\infty)\big)$ from $0$ to $z$. It extends continuously to $0$ and $1$ with
$
\mathcal{L}(0)=\Li_2(0)=0,$ and $ \mathcal{L}(1)=\Li_2(1)=\frac{\pi^2}{6}$. We have the following closed form of $I(a,b)$. 

\begin{lemma}
\label{lem:LassoCone2}
Let $0<a<b<1$. Then
$$
\begin{aligned}
I(a,b)
&=\mathcal{L}(z)+\mathcal{L}(\overline z)-\mathcal{L}(1/z)-\mathcal{L}(1/\overline{z})\\
&=2\,\mathrm{Re}\left\{\mathcal{L}(z)-\mathcal{L}(1/z)\right\},\quad \text{where } z:=\frac{\frac{1-ab}{a+b}+i}{\frac{1-b^2}{2b}+i}. 
\end{aligned}
 $$
\end{lemma}

\begin{proof}
Set $a=\tan\alpha$, $b=\tan\beta$ with $0<\alpha<\beta<\tfrac{\pi}{4}$, and substitute $x=\tan u$, $y=\tan v$.
Using $1+\tan^2 t=\sec^2 t$ and $\tan t-\tan s=\frac{\sin(t-s)}{\cos t\cos s}$ we obtain
\[
\frac{dx\,dy}{(y-x)^2}=\frac{du\,dv}{\sin^2(v-u)},\qquad
\frac{(x-a)(x-b)(1+y^2)}{(y-a)(y-b)(1+x^2)}
=\frac{\sin(u-\alpha)\sin(u-\beta)}{\sin(v-\alpha)\sin(v-\beta)}.
\]
Define
$$
\Phi(t):=\ln\big|\sin(t-\alpha)\big|+\ln\big|\sin(t-\beta)\big|
$$
where $\ln$ is the natural logarithm defined on reals.

Then
\begin{equation}\label{eq:I-angle}
I:=I(a,b) = \int_{\alpha}^{\beta}\int_{v-\pi/2}^{-\beta}
\frac{\Phi(u)-\Phi(v)}{\sin^2(v-u)}\,du\,dv.
\end{equation}
Integrating by parts with respect to $u$ using $\partial_u\big[\cot(v-u)\big]=1/\sin^2(v-u)$, we obtain 
$$
\int \frac{\Phi(u)-\Phi(v)}{\sin^2(v-u)}\,du
=\Big[(\Phi(u)-\Phi(v))\cot(v-u)\Big]
-\int \Phi'(u)\,\cot(v-u)\,du,
$$
where $\Phi'(t)=\cot(t-\alpha)+\cot(t-\beta)$.
At the lower endpoint $u=v-\frac{\pi}{2}$, $\cot(v-u)=\cot(\frac{\pi}{2})=0$.
Therefore
\begin{equation}\label{eq:after-IBP}
I=I_1-I_2,
\end{equation}
where $$I_1:=\int_{\alpha}^{\beta}\Big(\Phi(-\beta)-\Phi(v)\Big)\cot(v+\beta)\,dv,\quad
I_2:=\int_{\alpha}^{\beta}\int_{v-\pi/2}^{-\beta}\Phi'(u)\,\cot(v-u)\,du\,dv.$$

Splitting and changing variables in $I_1$, we obtain
\begin{align}
I_1&=\Phi(-\beta)\int_{\alpha}^{\beta}\cot(v+\beta)\,dv-\int_{\alpha}^{\beta}\Phi(v)\cot(v+\beta)\,dv\\
&=2F_0(\alpha+\beta,2\beta) -F_{\alpha+\beta}(0,\beta-\alpha)+F_{-2\beta}(0,\beta-\alpha), \label{eq:first-part}
\end{align}
where 
\begin{equation}\label{eq:Ffunction}
 F_{\delta}(r_1,r_2):=\int_{r_1}^{r_2}(\ln \sin t)\cot(t+\delta)\,dt
,
\end{equation}
and for $\delta=0$ one has the limit form
$F_{0}(r_{1},r_{2})=\frac12\bigl(\ln^{2}\sin r_{2}-\ln^{2}\sin r_{1}\bigr).
$

Now, consider the trapezoidal domain
$
\mathcal R=\{(u,v):\ \alpha\le v\le\beta,\ v-\tfrac{\pi}{2}\le u\le -\beta\},
$
we have
$$
v\in
\begin{cases}
[\,\alpha,\,u+\tfrac{\pi}{2}\,], & u\in[\alpha-\tfrac{\pi}{2},\,\beta-\tfrac{\pi}{2}],\\
[\,\alpha,\,\beta\,], & u\in[\,\beta-\tfrac{\pi}{2},\, -\beta\,].
\end{cases}
$$
Applying the Fubini theorem to swap the order of integration in $I_2$, and then splitting, merging, and changing variables, we obtain
\begin{align}
I_2
&=\int_{\alpha - \pi/2}^{\beta - \pi/2} \Phi'(u) \left[ \int_{\alpha}^{u+\pi/2} \cot(v-u)\, dv \right] du +
\int_{\beta - \pi/2}^{-\beta} \Phi'(u) \left[ \int_{\alpha}^{\beta} \cot(v-u) dv \right] du\\
&=-\int_{\alpha-\pi/2}^{-\beta}\Phi'(u)\,\ln\sin(\alpha-u)\,du
+\int_{\beta-\pi/2}^{-\beta}\Phi'(u)\,\ln\sin(\beta-u)\,du \\
&=F_0(\alpha+\beta,2\beta)
+F_{\beta-\alpha}(\alpha+\beta,\tfrac{\pi}{2})-F_{-\beta+\alpha}(2\beta,\tfrac{\pi}{2}),\label{eq:J-final}
\end{align}
where $F$ is defined in \cref{eq:Ffunction}. Combining \cref{eq:after-IBP,eq:J-final,eq:first-part}, we get
\begin{equation}\label{eq:IinF}
I(a,b)=F_0(\alpha+\beta,2\beta) -F_{\alpha+\beta}(0,\beta-\alpha)+F_{-2\beta}(0,\beta-\alpha)-F_{\beta-\alpha}(\alpha+\beta,\tfrac{\pi}{2})+F_{-\beta+\alpha}(2\beta,\tfrac{\pi}{2}). 
\end{equation}

\begin{claim}\label{claim:antideriv}
Let $0\le r_1<r_2\le \tfrac{\pi}{2}$ and $\delta\in(-\tfrac{\pi}{2},\tfrac{\pi}{2})$ with $-\delta\notin[r_1,r_2]$. Then
$$\begin{aligned}
F_\delta(r_1,r_2)
:&=\int_{r_1}^{r_2}(\ln\sin t)\cot(t+\delta)\,dt\\
&=-\tfrac12\sum_{\varepsilon=\pm1}\Big[\Li_2\Big(\tfrac{\cot\delta-\cot(r_2+\delta)}{\cot\delta+\varepsilon i}\Big)
-\Li_2\Big(\tfrac{\cot\delta-\cot(r_1+\delta)}{\cot\delta+\varepsilon i}\Big)\Big]+\Phi_{\delta}(r_2)-\Phi_{\delta}(r_1),
\end{aligned}$$
where $\qquad \Phi_{\delta}(r):=\ln\Big(\tfrac{|\sin(r+\delta)|}{|\sin\delta|}\Big)\ln(\sin r)
-\tfrac12\ln^2\Big(\tfrac{|\sin(r+\delta)|}{|\sin\delta|}\Big)$.
\end{claim}

\begin{proof}[Proof of \cref{claim:antideriv}]
The case $\delta=0$ is given in \cref{eq:Ffunction}, thus we assume $\delta\neq 0$.
Set $u:=\cot\delta$ and $v_j:=\cot(r_j+\delta)$ for $j=1,2$. With $x:=\cot(t+\delta)$ we have $dx=-(1+x^2)\,dt$ and
$$
\cot t=\frac{ux+1}{u-x},\qquad \sin^2 t=\frac{(u-x)^2}{(1+u^2)(1+x^2)}.
$$
Thus
$$
\ln(\sin t)=\ln|u-x|-\tfrac12\ln(1+u^2)-\tfrac12\ln(1+x^2),
$$
and
$$
\begin{aligned}
F_\delta(r_1,r_2)
&=-\int_{v_1}^{v_2}\frac{x}{1+x^2}\Bigl(\ln|u-x|-\tfrac12\ln(1+u^2)-\tfrac12\ln(1+x^2)\Bigr)\,dx\\
&=J(u;v_1,v_2)+\tfrac14\ln(1+u^2)\ln\frac{1+v_2^2}{1+v_1^2}
+\tfrac18\bigl(\ln^2(1+v_2^2)-\ln^2(1+v_1^2)\bigr),
\end{aligned}
$$
where
$$
J(u;v_1,v_2):=\int_{v_2}^{v_1}\frac{x}{1+x^2}\ln|u-x|\,dx
=\tfrac12\bigl(R(i)+R(-i)\bigr),\qquad
R(z):=\int_{v_2}^{v_1}\frac{\ln|u-x|}{x-z}\,dx.
$$
Define
$$
\Psi_z(x):=\Li_2\Bigl(\frac{u-x}{u-z}\Bigr)+\ln|u-x|\Log\Bigl(\frac{x-z}{u-z}\Bigr), \,\,\,   z=\pm i.$$
Since $u\in\mathbb R$ and $z=\pm i$, the paths $\eta(x)=\frac{x-z}{u-z}$ and $\xi(x)=\frac{u-x}{u-z}$ meet the real axis only at $\eta(u)=1$ and $\xi(u)=0$. Hence $\eta(x)\notin(-\infty,0]$ and $\xi(x)\notin[1,\infty)$ for all real $x$ in the interval, so the principal branches of $\Log$ and $\Li_2$ are valid along these paths. Therefore the chain rule with $\frac{d}{d\xi}\Li_2(\xi)=-\frac{\Log(1-\xi)}{\xi}$ applies, so a direct differentiation yields
$$
\Psi_z'(x)=\frac{\ln|u-x|}{x-z}.
$$
Consequently, $R(z)
=\Psi_z(v_1)-\Psi_z(v_2)$, and hence 
$$
\begin{aligned}
J(u;v_1,v_2)
&=-\tfrac12\sum_{\varepsilon=\pm1}\Bigl[\Li_2\Bigl(\frac{u-v_2}{u-\varepsilon i}\Bigr)-\Li_2\Bigl(\frac{u-v_1}{u-\varepsilon i}\Bigr)\Bigr]\\
&\quad-\tfrac12\Bigl[\ln|u-v_2|\ln\Bigl(\frac{1+v_2^2}{1+u^2}\Bigr)-\ln|u-v_1|\ln\Bigl(\frac{1+v_1^2}{1+u^2}\Bigr)\Bigr].
\end{aligned}
$$

Now use $1+u^2=\csc^2\delta$ and $1+v_j^2=\csc^2(r_j+\delta)$ and the identity
$$
\ln(\sin r_j)=\ln|u-v_j|-\tfrac12\ln(1+u^2)-\tfrac12\ln(1+v_j^2),
\qquad j=1,2,
$$
to rewrite the logarithmic contribution in $F_\delta(r_1,r_2)$ as
$$
\Phi_\delta(r_2)-\Phi_\delta(r_1),
\qquad
\Phi_{\delta}(r):=\ln\Bigl(\frac{|\sin(r+\delta)|}{|\sin\delta|}\Bigr)\ln(\sin r)
-\tfrac12\ln^2\Bigl(\frac{|\sin(r+\delta)|}{|\sin\delta|}\Bigr).
$$
If $r_j=0$, interpret $\Phi_\delta(0)$ by its continuous extension, namely $\Phi_\delta(0)=0$.

Finally, replace $u=\cot\delta$ and $v_j=\cot(r_j+\delta)$  to obtain the stated formula of $F_\delta(r_1,r_2)$.
\end{proof}

Applying \cref{claim:antideriv} to \cref{eq:IinF}, the logarithmic part $\Phi_{\delta}(r)$ cancels, and the expression reduces to the following eight $\Li_2$ terms:
\begin{equation}\label{eq:8terms}
I(a,b)=\tfrac12\sum_{\varepsilon=\pm1}\Big[
\Li_2\Big(1-\tfrac{1}{z_{\varepsilon}}\Big)
-\Li_2\Big(1-z_{\varepsilon}\Big)
-\Li_2\Big(1/z_{\varepsilon}\Big)
+\Li_2\Big(z_{\varepsilon}\Big)\Big], 
\end{equation}
where
$$
z_{\varepsilon}:=\frac{\cot(\beta+\alpha)+\varepsilon i}{\cot(2\beta) +\varepsilon i}=\tfrac{\frac{1-ab}{a+b}+\varepsilon i}{\frac{1-b^2}{2b}+\varepsilon i}\,\,, \qquad \text{for $\varepsilon\in\{\pm1\}$.}
$$
Applying Euler's reflection identity $\Li_2(1-z)=\tfrac{\pi^2}{6}-\Li_2(z)-\Log (z)\, \Log(1-z)$ together with the definition $\mathcal{L}(z):=\Li_2(z)+\tfrac12\Log (z)\,\Log(1-z)$ and the conjugation property $\mathcal{L}(z)=\overline{\mathcal{L}(\overline{z})}$ (valid on the principal branches of $\Log$ and $\Li_2$) yields the desired formula.
\end{proof}
Our \cref{pro:LassoCuspCone} is a combination of \cref{lem:LassoCusp,lem:LassoCone1,lem:LassoCone2}.

\section{Full Liouville measure of generic vectors}\label{sec:Nongeneric}
Let $X$ be a hyperbolic cone surface obtained from a complete hyperbolic cone surface $\overline X$ by deleting all the cone points. Let $A_{\mathrm{cones}}\;:=\;\overline X\setminus X$ denote the set of cone points with positive angles. Let $\widetilde X$ be the universal cover of $X$, with covering projection 
\[
\pi^{\circ}\colon\widetilde X\longrightarrow X,
\]
endowed with the lifted hyperbolic metric, and let $\overline{\widetilde X}$ be the metric completion of $\widetilde X$. The map $\pi^{\circ}$ extends continuously to
\[
\pi\colon \overline{\widetilde X} \longrightarrow \overline X.
\]
We define the set of lifted cone points as $\widetilde A_{\mathrm{cones}} := \overline{\widetilde X} \setminus \widetilde X = \pi^{-1}(A_{\mathrm{cones}})$. \begin{remark}\label{rem:completion-cone-preimage} The metric completion $\overline{\widetilde X}$ is a simply connected space containing the discrete subset $\widetilde A_{\mathrm{cones}}$. The geometry near a lifted cone point $\tilde p \in \widetilde A_{\mathrm{cones}}$ is that of a hyperbolic cone with infinite total angle: metric circles centered at $\tilde p$ have infinite length. Consequently, $\overline{\widetilde X}$ is distinct from $\widetilde{\overline X}$, the universal orbifold cover of $\overline X$ (defined when $X$ is a good orbifold \cite[Chapter 13]{thurston2022geometry}). In $\widetilde{\overline X}$, the preimage of a cone point is a regular point with total angle $2\pi$, reflecting the finite order of the local orbifold fundamental group, whereas in $\overline{\widetilde X}$, the "logarithmic" unwrapping of the puncture creates a singularity of infinite angle. \end{remark}

Fix a developing map
$$
\mathrm{dev}^\circ\colon \widetilde X \longrightarrow \mathbb D
$$
for the hyperbolic metric on $X$. This map is $\pi_1(X)$–equivariant and a
local isometry. Since $\mathrm{dev}^\circ$ is $1$–Lipschitz, it extends
uniquely to a continuous $\pi_1(X)$–equivariant map
\begin{equation}\label{eq:developingMap}
\mathrm{dev}\colon \overline{\widetilde X} \longrightarrow \mathbb D,
\end{equation}
which is still a local isometry on
$\widetilde X = \overline{\widetilde X} \setminus \widetilde A_{\mathrm{cones}}$.

\begin{definition}[Endpoint map]\label{dfn:EndpointMap}
For $x\in \overline X$, choose a lift $\tilde x\in \overline{\widetilde X}$ and
set $\hat x:=\mathrm{dev}(\tilde x)\in\mathbb D$. Let $\Sigma_x \overline X$
denote the space of directions at $x$ (so $\Sigma_x \overline X\simeq T^1_x X$
if $x\in X$). Given $u\in\Sigma_x \overline X$, choose $\varepsilon>0$ and let
\[
\gamma_{x,u}\colon[0,\varepsilon)\to \overline X
\]
be the geodesic starting at $x$ with initial direction $u$ such that
$\gamma_{x,u}\bigl((0,\varepsilon)\bigr)\subset X$. Since
$\pi^\circ\colon \widetilde X\to X$ is a covering, there is a unique lift
\[
\tilde\gamma_{x,u}\colon (0,\varepsilon)\to \widetilde X
\]
satisfying $\pi^\circ\circ\tilde\gamma_{x,u}=\gamma_{x,u}|_{(0,\varepsilon)}$
and $\lim_{t\to 0}\tilde\gamma_{x,u}(t) = \tilde x$; we set
$\tilde\gamma_{x,u}(0):=\tilde x$. Along $(0,\varepsilon)$ the map $\mathrm{dev}$
is a local isometry, hence $\mathrm{dev}\circ\tilde\gamma_{x,u}$ is a geodesic
segment in $\mathbb D$ which extends uniquely to a geodesic ray
\[
r_{x,u}\colon[0,\infty)\to\mathbb D
\]
with $r_{x,u}(0)=\hat x$ and $r_{x,u}|_{(0,\varepsilon)}=
\mathrm{dev}\circ\tilde\gamma_{x,u}$. The endpoint map
\begin{equation}\label{eq:EndingMap}
\operatorname{end}_x\colon\Sigma_x \overline X\to\partial\mathbb D,\qquad
\operatorname{end}_x(u):=\lim_{t\to\infty} r_{x,u}(t)
\end{equation}
is then well defined for the chosen lift $\tilde x$. Changing $\tilde x$ by a
deck transformation $g\in\pi_1(X)$ post-composes $\operatorname{end}_x$ with the boundary action of the corresponding isometry of $\mathbb D$.
\end{definition}

Fix $\varepsilon_0>0$ sufficiently small so that the radius-$\varepsilon_0$ collars in $\overline X$ around all non-orbifold cone points are pairwise disjoint. Let $X_{\varepsilon_0}$ be the subsurface obtained from $\overline X$ by removing the union of these $\varepsilon_0$-collars of non-orbifold cone points. Then $X_{\varepsilon_0}$ is a concave core of $\overline X$.

For $p\in \overline X$, a direction $u\in\Sigma_p \overline{X}$ is called \emph{non-generic} with respect to the concave core $X_{\varepsilon_0}$ if it satisfies \cref{dfn:NongenericVector}; otherwise, it is \emph{generic}. When $p$ is smooth, non-genericity is antipodally invariant: for any $v\in T_p^1 X$, $v$ is non-generic $\iff$ $-v$ is non-generic. The following theorem implies that the set of generic directions has full Liouville measure.

\begin{theorem}\label{pro:WeakErgodicity}
Let $X_{\varepsilon_0}$ be as defined above, and for any $p\in \overline X$ let $N_p$ denote the set of non-generic directions at $p$ with respect to the concave core $X_{\varepsilon_0}$. 
If $\theta$ denotes the angular (Lebesgue) measure on $\Sigma_p \overline{X}$ (the space of directions at $p$), then $\theta(N_p)=0$.
\end{theorem}
\begin{proof} If $p$ is a cone point or lies on a geodesic boundary component, the claim follows directly from \cite[Theorem~5.3]{basmajian2025prime}.

Assume, for contradiction, that $\theta(N_p)>0$ for some smooth point $p\in X$.
Identify $N_p\subset T^1_p(X)$ with the set of geodesic rays from $p$ determined by these vectors.
Since $N_p$ is antipodally symmetric (i.e.\ $v\in N_p \Leftrightarrow -v\in N_p$) and $\theta(N_p)>0$, there exists a sector $\delta$ at $p$ of arbitrarily small angle $|\delta|$ with $\theta(N_p\cap\delta)>0$.
If $p$ lies in the $\varepsilon_0$–collar of a non-orbifold cone point $c$, choose $\delta$ so that every ray in $\delta$ initially increases its distance from $c$.
Fix $0<\varepsilon_1<\varepsilon_0$ so that
\[
d_X(p,c)>\varepsilon_1 \quad \text{for every } c\in A_{\mathrm{non-orb}},\]
where $A_{\mathrm{non-orb}}$ denotes the set of all non-orbifold cone points on $\overline X$.
Define $X_{\varepsilon_1}$ analogously to $X_{\varepsilon_0}$ by removing the union of radius-$\varepsilon_1$ collars around all non-orbifold cone points.
By the choice of $\delta$ and \cref{dfn:NongenericVector}, all geodesic rays in $N_p\cap\delta$ lie entirely inside $X_{\varepsilon_1}$.
By compactness (with $|\delta|$ small enough), we may moreover assume that the two boundary geodesic rays $\alpha_1,\alpha_2$ of $\delta$ (whose angular separation equals the aperture of $\delta$) also lie entirely inside $X_{\varepsilon_1}$.

\begin{figure}[h!]
 \centering
 \begin{tikzpicture}[scale=0.4]


 \draw[thin] (0,0) circle (8);

 \draw[thin,dashed] (-2,0) -- (2,0);

 \filldraw[fill=gray!20,draw=black,thin]
 (0,0) -- (-75:8)
 arc[start angle=-75, end angle=-105, radius=8]
 -- cycle;

 \draw[thin] (0,0) -- (-75:8);
 \draw[thin] (0,0) -- (-105:8);

 \draw[thin] (-105:2) arc[start angle=-105, end angle=-75, radius=2];
 \node at (-90:2.5) {$\hat \delta$};

 \draw[thin] (0,-0.5) -- (0,0.5);

 \draw[thin, dashed] (0,0) -- (1.5,2);

 \begin{scope}[shift={(-17,-3.726)}]
 \draw[thin, dashed] (15,3.726)
 arc[start angle=14, end angle=-15, radius=15.458];
 \end{scope}

 \begin{scope}[shift={(17,-3.726)}]
 \draw[thin, dashed] (-15,3.726)
 arc[start angle=166, end angle=195, radius=15.458];
 \end{scope}

 \begin{scope}[shift={(-17,0)}]
 \draw[thin, dashed] (15,0)
 arc[start angle=0, end angle=28.07, radius=15];
 \end{scope}
 
 \begin{scope}[shift={(17,0)}]
 \draw[thin, dashed] ({15*cos(151.93)},{15*sin(151.93)})
 arc[start angle=151.93, end angle=180, radius=15];
 \end{scope}

 \draw[thin] (-2,0) -- ++(0,0.4) -- ++(0.4,0) -- ++(0,-0.4);
 \draw[thin] (2,0) -- ++(0,0.4) -- ++(-0.4,0) -- ++(0,-0.4);

 \begin{scope}[shift={(0,17.5625)}]
 \draw[thin, dashed]
 ({15.634622*cos(-97.883)},{15.634622*sin(-97.883)})
 arc[start angle=-97.883, end angle=-82.117, radius=15.634622];
 \end{scope}

 \begin{scope}[shift={(25.38899744,-1.47924808)}]
 \draw[thin]
 ({24.14103076*cos(171.71357385)},{24.14103076*sin(171.71357385)})
 arc[start angle=171.71357385, end angle=195, radius=24.14103076];
 \end{scope}
 \begin{scope}[shift={(53.61407768,-22.64805826)}]
 \draw[thin]
 ({57.64897110*cos(154.68755823)},{57.64897110*sin(154.68755823)})
 arc[start angle=154.68755823, end angle=165, radius=57.64897110];
 \end{scope}

 \begin{scope}[shift={(17,4.8125)}]
 \draw[thin, dashed]
 ({15.75309989*cos(-169.7154856)},{15.75309989*sin(-169.7154856)})
 arc[start angle=-169.7154856, end angle=-162.2120275, radius=15.75309989];
 \end{scope}

 \begin{scope}[shift={(-17,30.3125)}]
 \draw[thin, dashed]
 ({33.82081691*cos(-63.67170689)},{33.82081691*sin(-63.67170689)})
 arc[start angle=-63.67170689, end angle=-56.83853051, radius=33.82081691];
 \end{scope}

 \fill (0,0) circle (2.5pt) node[above] {$\widehat p$};
 \fill (-2,0) circle (2.5pt) node[left] {$\widehat p_1$};
 \fill (2,0) circle (2.5pt) node[right] {$\widehat p_2$};
 \fill (-75:8) circle (2.5pt) node[above left] {$\widehat \alpha_2$};
 \fill (-105:8) circle (2.5pt) node[above right] {$\widehat \alpha_1$};
 \node at (-2,-8.5) {$e_1$};
 \node at (2,-8.5) {$e_2$};
 \filldraw[fill=white,draw=black] (1.5,2) circle (2.5pt) node[above] {$\widehat q$};

 \node at (-2.5,-3) {$\widehat \alpha'_1$};

 \node at (2.5,-3) {$\widehat \alpha'_2$};

 \end{tikzpicture}
\caption{Images in  $\mathbb D$ of relevant lifts under the developing map.}

 \label{fig:Developing}
\end{figure}
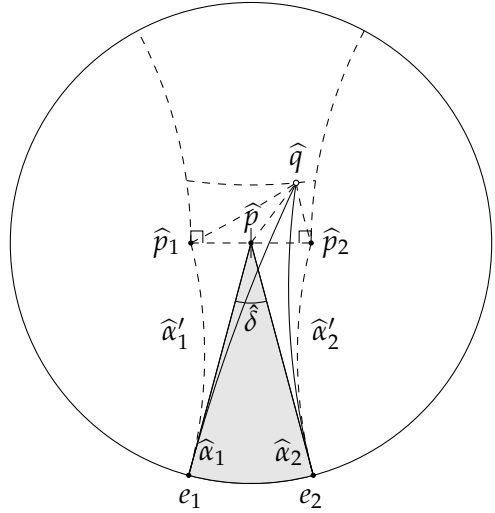


Let \(\beta:[-\varepsilon_1,\varepsilon_1]\to X\) be the unit-speed oriented geodesic segment with midpoint \(\beta(0)=p\), lying between $\alpha_1$ and $\alpha_2$, which bisects $\delta$. Let $\eta$ be the geodesic segment of length $\tfrac{\varepsilon_1}{2}$ perpendicular to $\beta$ at $p$, with endpoints $p_1,p_2$ and midpoint $p$.
Transport $\eta$ at unit speed along the geodesic flow orthogonal to $\eta$ in the negative
$\beta$–direction until the first intersection time $\tau$ with
a boundary element of $X$; choose a first hitting point $q\in\eta_\tau$. The translates $\{\eta_s\}_{0\le s\le\tau}$ sweep out a
compact immersed strip $K$ whose only possible contacts with the set of boundary elements of $X$ occur on the terminal translate $\eta_\tau$. Let $[q,p]$, $[q,p_1]$, and $[q,p_2]$ be the oriented geodesic segments embedded in $K$
from $q$ to $p$, $p_1$, and $p_2$, respectively.

Since $\alpha_i\subset X_{\varepsilon_1}$ and the length of the geodesic segment $[p,p_i]$ is $\frac{\varepsilon_1}{4}$ ($i\in\{1,2\}$),
the concatenation $[p_i,p]\ast\alpha_i$ is contained in $X_{\frac{3}{4}\varepsilon_1}$ and is homotopic (in the orbifold sense) to a geodesic ray issuing from $p_i$, which we denote by $\alpha_i'$. 

For $|\delta|$ sufficiently small, we may assume that, for each $i\in\{1,2\}$, the interior angle at $p_i$ satisfies
\[
\measuredangle_{p_i}\bigl([p_i,p],\alpha_i'\bigr)\ \le\ \frac{\pi}{2}.
\]
Similarly, for each $\eta\in N_p\cap\delta$ and each $i\in\{1,2\}$, the concatenation
$[p_i,p]\ast\eta$
is homotopic (in the orbifold sense) to a geodesic ray issuing from $p_i$, which we denote by $\eta_i'$. In particular, the rays $\eta$ and $\eta_i'$ are asymptotic, and the associated geodesic triangle with one ideal vertex, immersed away from orbifold cone points, lies entirely inside $X_{\frac{3}{4}\varepsilon_1}$.

By the above construction and compactness, there exists
$0<\varepsilon_2<\varepsilon_1$ such that, for every $\eta\in N_p\cap\delta$, the convex quadrilateral, immersed away from orbifold cone points, bounded by $[q,p_1]$, $[q,p_2]$, $\eta_1'$, and $\eta_2'$ lies in $X_{\varepsilon_2}$, except possibly within the $\varepsilon_2$–collar about $q$ when $q\in A_{\mathrm{non-orb}}$.
 As a consequence of convexity, there exists uniquely a
geodesic ray $\gamma$ with $\gamma(0)=q$ in the homotopy class (in the orbifold sense) of the concatenation $[q,p]\ast\eta$, and $\gamma$ remains inside $X_{\varepsilon_2}$ after its first
exit from the $\varepsilon_2$-collar centered at $q$ (if $q \in A_{\mathrm{non-orb}}
$). Define $N_{p\to q}$ to be the set of such geodesic rays $\gamma$. Identify $N_{p\to q}$ with the corresponding subset of $\Sigma_q \overline{X}$ (the space of directions at $q$). Let $\delta'\subset\Sigma_q \overline{X}$ be the minimal sector containing $N_{p\to q}$. Then
$N_{p\to q}\subset N_q\cap\delta'$,
where $N_q\subset\Sigma_q \overline{X}$ denotes the non-generic directions at $q$ (with respect to $X_{\varepsilon_2}$). Moreover, for $|\delta|$ small enough, $N_{p\to q}$ is in bijection with $N_p\cap\delta$ via
\begin{equation}\label{eq:Phi}
\Phi_{p,q}:N_p\cap\delta\to N_{p\to q},\qquad 
\Phi_{p,q}(\eta)=\gamma\ \text{ where }\ [q,p]*\eta \simeq \gamma
\text{ in the orbifold sense}.\end{equation}

Postcomposing the developing map $\mathrm{dev}$ (see \cref{eq:developingMap}) with an isometry of $\mathbb D$, we choose a lift
$\tilde p\in\widetilde X$ of $p$ such that $\hat p:=\mathrm{dev}(\tilde p)=0$. See \cref{fig:Developing}
for the images under $\mathrm{dev}$ of the relevant lifts of $p$, $q$, $p_1$, $p_2$, as well as the sectors $\delta$ and $\delta'$
constructed above, together with the corresponding geodesic rays.
 For these sectors, we have
\[
\operatorname{end}_p(\delta)=\operatorname{end}_q(\delta')=:I\subset\partial\mathbb D,
\]
where the endpoint map $\operatorname{end}_x(u)$ is defined in \cref{dfn:EndpointMap}, and hence
\[
\Phi_{p,q}:=\bigl(\operatorname{end}_q\vert_{\delta'}\bigr)^{-1}\circ\bigl(\operatorname{end}_p\vert_{\delta}\bigr)
\colon \delta\xrightarrow{\ \sim\ }\delta'
\]
is a bijection with $\operatorname{end}_q\circ\Phi_{p,q}=\operatorname{end}_p$ on $\delta$.
On $N_p\cap\delta$, this coincides with the homotopy definition in \cref{eq:Phi}. We have the following standard fact in hyperbolic geometry.
\begin{lemma}\label{lem:biLip-local-0}
The map \(\Phi_{p,q}:\delta\to\delta'\) is bi–Lipschitz for the angular metrics at \(p\) and \(q\). In particular,
\[
\frac{1}{L}\,\measuredangle_p(u,v)\ \le\ \measuredangle_q\bigl(\Phi_{p,q}(u),\Phi_{p,q}(v)\bigr)\ \le\ L\,\measuredangle_p(u,v)\qquad(u,v\in\delta),
\]
where \(L:=\tfrac{\pi}{2}\,e^{\ell_X([p,q])}\).
\end{lemma}

\begin{proof}[Proof of \cref{lem:biLip-local-0}]
Let $\ell:=\ell_X([p,q]) = d_{\mathbb D}(0,\hat q)$. Write \(M_z(w)=\dfrac{z+w}{1+\overline z\,w}\) and let \(d_0(\xi,\eta)=|\xi-\eta|\) be the chordal metric on
\(\partial\mathbb D\).
Then
\[
d_p(\xi,\eta):=d_0(\xi,\eta),\qquad
d_q(\xi,\eta):=d_0 \bigl(M_{\hat q}^{-1}\xi,M_{\hat q}^{-1}\eta\bigr),
\]
and the boundary derivative estimate
\[
\Bigl|(M_{\hat q}^{-1})'(e^{i\theta})\Bigr|\in\Bigl[\frac{1-|\hat q|}{1+|\hat q|},\,\frac{1+|\hat q|}{1-|\hat q|}\Bigr]
=\bigl[e^{-\ell},\,e^{\ell}\bigr]
\]
gives
\begin{equation}\label{eq:distort}
e^{-\ell}\,d_p(\xi,\eta)\ \le\ d_q(\xi,\eta)\ \le\ e^{\ell}\,d_p(\xi,\eta).
\end{equation}
For \(\xi,\eta\in\partial\mathbb D\) with central angle \(\measuredangle_0(\xi,\eta)\in[0,\pi]\),
\[
\frac{2}{\pi}\,\measuredangle_0(\xi,\eta)\ \le\ d_p(\xi,\eta)=2\sin \Bigl(\frac{\measuredangle_0(\xi,\eta)}{2}\Bigr)\ \le\ \measuredangle_0(\xi,\eta),
\]
and likewise \(\measuredangle_q(\xi,\eta)\le \tfrac{\pi}{2}\,d_q(\xi,\eta)\).
Now for \(u,v\in\delta\),
\[
\measuredangle_q \bigl(\Phi_{p,q}(u),\Phi_{p,q}(v)\bigr)
\ \le\ \tfrac{\pi}{2}\,d_q \bigl(\operatorname{end}_p u,\operatorname{end}_p v\bigr)
\ \le\ \tfrac{\pi}{2}\,e^{\ell}\,d_p \bigl(\operatorname{end}_p u,\operatorname{end}_p v\bigr)
\ \le\ \tfrac{\pi}{2}\,e^{\ell}\,\measuredangle_p(u,v),
\]
using \cref{eq:distort} and \(d_p\le\measuredangle_p\).
The reverse inequality follows by swapping \(p\) and \(q\).
\end{proof}
In particular, for any measurable \(E\subset\delta\),
\(\theta(E)=0\ \Longleftrightarrow\ \theta\bigl(\Phi_{p,q}(E)\bigr)=0\). By \cite[Theorem~5.3]{basmajian2025prime} and the fact $N_{p\to q}=\Phi_{p,q}(N_p\cap\delta)\subset N_q$, we have that $$ \theta(N_q)=0
\ \Rightarrow\ \theta(N_{p\to q})=0
\ \Rightarrow\ \theta(N_p\cap\delta)=0,$$ which contradicts the assumption $\theta(N_p\cap\delta)>0$.
\end{proof}
\begin{remark}
Let $V:=V_{\vec{k}}(X)$ denote the concave core of $X$, where $(X,\vec{k})$ is an admissible pair. By \cref{dfn:kNaturalCollar}, we can always choose $\varepsilon_0 >0$ sufficiently small so that $V\subseteq X_{\varepsilon_0}$. Consequently, the set of non-generic vectors with respect to $V$ is also of measure zero.
\end{remark}

{\small
\bibliographystyle{alpha}
\bibliography{references}
}

The Graduate Center, CUNY, 365 Fifth Ave., New York, NY 10016, {\it and}\\
Hunter College, CUNY, 695 Park Ave., New York, NY 10065, USA\\
{\em Email address:}
 \href{mailto:abasmajian@gc.cuny.edu}{abasmajian@gc.cuny.edu}
 
Department of Mathematics, National University of Singapore, Singapore, {\it and}\\
Institute of Mathematics, Vietnam Academy of Science and Technology, Vietnam\\
{\em Email address:}
 \href{mailto:dnminh@math.ac.vn}{dnminh@math.ac.vn}
 
Department of Mathematics, University of Fribourg, Switzerland\\
{\em Email address:}
\href{mailto:hugo.parlier@unifr.ch}{hugo.parlier@unifr.ch}

Department of Mathematics, National University of Singapore, Singapore\\
{\em Email address:}
\href{mailto:mattansp@nus.edu.sg}{mattansp@nus.edu.sg}

\end{document}